\documentclass{amsproc}

\usepackage{graphicx}

\newtheorem{The}{Theorem}[section]
\newtheorem{Lem}[The]{Lemma}
\newtheorem{Pro}[The]{Proposition}
\newtheorem{Cor}[The]{Corollary}

\theoremstyle{definition}
\newtheorem{Def}[The]{Definition}

\theoremstyle{remark}
\newtheorem{Rem}[The]{Remark}

\numberwithin{equation}{section}

\newcommand{\real}{ {\mathbb R} }
\newcommand{\complex}{ {\mathbb C} }
\newcommand{\integer}{ {\mathbb Z} }

\newcommand{\rational}{ {\mathbb Q} }

\newcommand{\fol}{ {\mathcal F} }

\newcommand{\fra}{ {\mathcal F} }
\newcommand{\hyp}{ {\mathbb H} }

\newcommand{\tw}{ {\tt t} }
\newcommand{\tb}{ {\tt tb} }
\newcommand{\rot}{ {\tt r} }

\begin{document}

\title{On the classification of tight contact structures}

\author{Paolo Ghiggini}
\address{Scuola Normale Superiore, Pisa}
\email{ghiggini@mail.dm.unipi.it}

\author{Stephan Sch{\"o}nenberger}
\address{University of Pennsylvania, Philadelphia}
\email{stephans@math.upenn.edu}

\date{January 11, 2002 and, in revised form, July 29, 2002.}



\begin{abstract}
Recently, there have been several breakthroughs in the classification of tight contact structures. We give an outline on how to exploit methods developed by Ko Honda and John Etnyre to obtain classification results for specific examples of small Seifert manifolds.
\end{abstract}

\maketitle

\section{Introduction}

After Eliashberg proved a classification for so-called overtwisted contact structures \cite{eliashberg:4},
work concentrated on the classification of tight contact structures, which turned out
to be much more subtle and provide interesting relations to the topology of the underlying manifold.
See \cite{etnyre:1} for an introduction to contact geometry and further references.

Until recently, the main tool to show that a contact structure on a manifold is tight is to show it is fillable.
A contact structure is holomorphically fillable if it is the oriented boundary of a compact
Stein 4-manifold. Gromov and Eliashberg showed that a fillable contact structure is
tight \cite{gromov, eliashberg:2}. Moreover, fillability is preserved by Legendrian surgery
\cite{weinstein, eliashberg:3}, thus providing a rich source of tight contact structures.
Gompf's extensive study on Legendrian surgery \cite{gompf:1} enables one in particular to construct holomorphically
fillable contact structures on many Seifert manifolds.
Using Legendrian surgery and techniques from Seiberg-Witten-theory, Lisca and Mati{\'c} \cite{lisca_matic}
proved that for every integer $n > 1$ there exist at least $[\frac{n}{2}]$ tight contact structures on the Brieskorn
homology spheres with reversed orientation, $-\Sigma(2,3,6n-1) = M(-\frac{1}{2}, \frac{1}{3}, \frac{n}{6n-1})$. Later, 
they improved this lower bound to $n-1$ in \cite{lisca_matic:2}.

Contact structures induce a singular foliation on embedded surfaces and these are often
easier to study than the contact structure itself. Motivated by work of
Eliashberg and Gromov \cite{eliashberg_gromov},  Giroux introduced the notion of convex surfaces, i.e.
~surfaces whose characteristic foliation is cut transversely by a certain multicurve, called
the dividing set \cite{giroux:5}. This dividing set essentially determines the contact structure in a
neighbourhood of the surface and is a convenient tool to study contact structures.

Exploiting this idea, Kanda gave a complete classification of tight contact structures on the 3-torus
\cite{kanda:1}; see also \cite{giroux:3}. 
This led Honda to study so-called bypasses attached along convex surfaces, which provide a systematic
tool for altering the dividing set of a convex surface. By splitting a contact 3-manifold along convex
surfaces into simpler pieces and studying the possibilities of tight contact structures, Honda
gave a complete classification of tight contact structures on solid tori, toric annuli, Lens spaces in 
\cite{honda:1}, as well as torus bundles over the circle and circle bundles over closed Riemannian surfaces;
see \cite{honda:2}. Many of these were independently obtained by Giroux \cite{giroux:3, giroux:2, giroux:1},
and on some Lens spaces by Etnyre \cite{etnyre:4}.

Furthermore, Lisca proved in \cite{lisca:3} that the Poincar{\'e} homology sphere with reverse orientation $-\Sigma(2,3,5)$ 
(this corresponds to the Seifert manifold $M(-\frac{1}{2}, \frac{1}{3}, \frac{1}{5})$ in the notation below)
has no symplectically (weakly) semi-fillable contact structure, thus proving  a conjecture of Gompf in
\cite{gompf:1}. Using the bypass technique in contact topology, Etnyre and Honda finally proved
\cite{etnyre_honda:1} the nonexistence of a tight contact structure on
$M(-\frac{1}{2},\frac{1}{3}, \frac{1}{5})$, thereby providing the first example of a closed 3-manifold which
admits no tight contact structure.

Lisca \cite{Lisca:2} went further and proved (among other things) that the Seifert manifolds
$M(-\frac{1}{2},\frac{1}{3},\frac{1}{4})$, $M(-\frac{1}{2},\frac{1}{3},\frac{1}{3})$ admit no (weakly) symplectically semi-fillable contact structure.
From Lisca's examples, Etnyre and Honda proved that on the Seifert manifolds
$M(-\frac{1}{2}, \frac{1}{4}, \frac{1}{4})$ and $M(-\frac{2}{3}, \frac{1}{3}, \frac{1}{3})$ there exist tight
contact structures without symplectic fillings.
These examples belong to the handful of Seifert manifolds which can be defined as torus bundles
over the circle; see \cite{geiges_ding:1} for fillability results on these manifolds.

Furthermore, the examples above are Seifert manifolds over the sphere $S^2$ with three singular fibres. On
`larger' Seifert manifolds it is recently proven by Colin \cite{colin:1}
that  every orientable Seifert manifold over a surface of genus $g\geq 1$ has infinitely many non-isomorphic
tight contact structures. Moreover, Colin \cite{colin:2}, see also Honda, Kazez, Mati{\'c} \cite{honda_kazez_matic:1}, proved that every
closed irreducible orientable toroidal 3-manifold carries infinitely many contact structures.

Therefore the classification of tight contact structures on Seifert manifolds may provide interesting 
new insight to the topology of tight contact structures on 3-manifolds. In this work, we will demonstrate
on two examples how to apply bypass techniques to obtain upper bounds on the number of tight contact structures
on Seifert manifolds over the sphere with three singular fibres. In the examples below, tight contact
structures are constructed using Legendrian surgery.

\section{Basic contact geometry}  

A positive contact structure on an oriented 3-manifold $M$ is a 2-plane field $\xi=\ker\alpha\subset TM$,
defined by a 1-form $\alpha$ satisfying $\alpha\land d\alpha > 0$. According to this definition, $\xi$
is co-oriented by $\alpha$ and oriented by $d\alpha$ such that the orientation on $M$ coincides with the orientation
defined by $\alpha\land d\alpha$.

\subsection{Legendrian curves and twisting}

A curve $\gamma$ in a contact manifold $(M,\xi)$ everywhere tangent to $\xi$ is called {\it Legendrian}.
Throughout this paper, we assume curves to be closed, and we will refer to `arcs' otherwise. Recall that every diffeomorphism between Legendrian curves extends to a contactomorphism of their neighbourhoods.
A Legendrian curve $\gamma$ in a contact manifold $(M, \xi)$ is endowed with a natural framing defined by a vector field along $\gamma$ transverse to $\xi$, called the {\it contact framing}. 
The {\it twisting number} $\tw(\gamma, \fra)$ is defined as the number of right $2\pi$ twists of the contact framing with respect to a preassigned framing $\fra$ of $\gamma$. 
In case $\gamma$ is a Legendrian boundary component of an oriented surface $S$, let $\fra_S$
denote the framing of $\gamma$ defined by $S$. In this case we will write $\tw(\gamma) := \tw(\gamma, \fra_S)$.
When $\gamma$ is null-homologous and $S$ is a Seifert surface for $\gamma$, the twisting number is called
{\it Thurston-Bennequin invariant} and denoted by $\tb(\gamma)$.

We have another classical invariant for Legendrian knots, the {\it rotation number} $\rot$:
If $\gamma$ is the Legendrian boundary of a Seifert surface $S$, we define $\rot(\gamma)$ as the number of revolution
of its tangent $\dot{\gamma}$ with respect to a trivialization of $\xi|_S$.
Note that, for any relative homology class $\beta \in H_2(M,\gamma)$ and $S\in\beta$ a representing surface, the
rotation number of $\gamma$ is independent of the choice of a trivialization but depends on $\beta$, and
reversing the orientation of $\gamma$ reverses the sign of $\rot$.
We refer the reader to \cite{etnyre:1} for a detailed discussion of the following
\begin{Pro}[Bennequin's Inequality]
        If $\gamma$ is a Legendrian knot in a tight contact manifold $(M, \xi)$ and $S$ a Seifert surface for $\gamma$ with Euler characteristic $\chi(S)$,
        then $$\tb(\gamma) +|\rot(\gamma)| \leq -\chi(S).$$ 
\end{Pro} 

\subsection{Convex surface theory}

Assume $S$ is a compact oriented surface embedded in a contact manifold $(M, \xi)$. The line field
$l_x = \xi_x \cap T_xS$, $x \in S$, integrates to a singular foliation $S_\xi$ of $S$ called
{\it characteristic foliation}. Recall that the singularities of $S_\xi$ are exactly the points in $S$ where
the contact plane is tangent to $S$. The characteristic foliation determines the contact structure in a
tubular neighbourhood and one has a certain freedom to alter the characteristic foliation by perturbing the
surface; see \cite{etnyre:1}. Generically, the amount of information needed to locally determine the contact
structure can be reduced to a collection of curves on the surface $S$.

A properly embedded orientable surface $S$ in a contact manifold $(M, \xi)$ is called {\it convex}, if
there exists a collection of curves $\Gamma$ on $S$ satisfying the following conditions:
\begin{enumerate}
        \item[(1)] $S\setminus\Gamma = S^+ \sqcup S^-$
        \item[(2)] $\Gamma$ is transverse to the characteristic foliation $S_\xi$ of $S$ 
        \item[(3)] There exists a vector field $v$ and a volume form $\theta$ on $S$ such that the
        characteristic foliation is directed by $v$, the flow of $v$ expands $\theta$ on $S^+$, contracts
        $\theta$ on $S^-$ and $v$ points transversely out of $S^+$. 
\end{enumerate}
Recall that the existence of {\it dividing curves} $\Gamma$ is equivalent to the existence of a contact vector
field $v$ transverse to the surface $S$, determining the contact structure in a neighbourhood of the surface up
to admissible isotopy, i.e. an isotopy $\phi:S\times [0,1] \to M$ such that $\phi(S\times \{t\})$ is transverse
to $v$ for all $t\in [0,1]$.

\begin{figure}[tb]
\includegraphics[width=7cm]{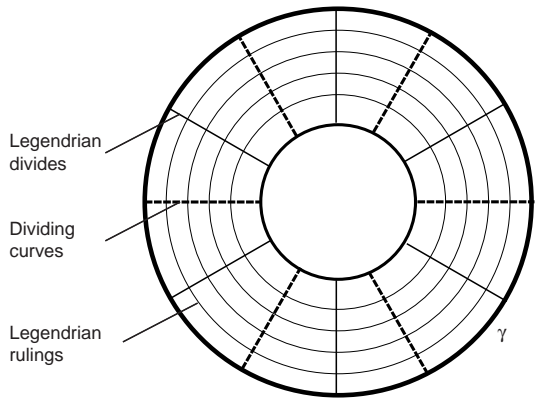}
\caption{Convex collar neighbourhood of a Legendrian boundary component with negative twisting.}
\label{f:collar}
\end{figure}

In \cite{giroux:5}, Giroux proved that every closed surface can be perturbed by a $C^\infty$-small
isotopy to be convex. More generally, a compact surface with Legendrian boundary can be perturbed to be convex
provided the twisting number of each boundary component is not positive. Moreover, the twisting number of a
boundary component $\partial S$ of a convex surface  $S$ determines the dividing set in a tubular
neighbourhood of $\partial S$. This follows from a relative version of Gray's Theorem in dimension three; see
Theorem 3.7 in \cite{etnyre:1}. We describe a standard tubular neighbourhood of a Legendrian boundary component
$\gamma$ as follows:
After perturbing $S$ we find a neighbourhood $N$ of a boundary component $\gamma \subset \partial S$ so that 
a collar neighbourhood $A=N\cap S$ of $\gamma$ in $S$ has the form
$A = S^1\times [0,1] =(\real/\integer)\times[0,1]$ with coordinates $(x,y)$ where  $\gamma= S^1\times\{0\}$.
In a neighbourhood $A\times[-1,1]$ of $A$  with  coordinates $(x,y,z)$  the
contact 1-form is defined by  $\alpha= \sin(2\pi n x)dy + \cos(2\pi n x)dz$ for $n = |\tw(\gamma)| \in \integer^+$.
Note that on this annulus the characteristic foliation consists of circles parallel to $\gamma$, called
{\it Legendrian rulings} and the dividing set consists of arcs transverse to the boundary, leading from one
boundary component to another. Between two dividing arcs lies an arc of  singularities, which we call
{\it Legendrian divides}; see Fig. \ref{f:collar}.
If $\tw(\gamma)=0$ then
the contact structure is defined by by $\alpha= dz-ydx$. In particular the twisting number of $\gamma$ is related to
the number of intersections of the dividing set $\Gamma$ with $\gamma$.

\begin{Pro} Suppose $S$ is a convex surface with Legendrian boundary in a contact manifold $(M, \xi)$
and $\gamma \in \partial S$ is a boundary component of $S$. Then
\begin{equation}\label{e:twcs}
\tw(\gamma) = -\frac{1}{2}\#(\gamma\cap\Gamma),
\end{equation}
where  $\#(\gamma\cap\Gamma)$ denotes the cardinality of the intersection $\gamma\cap\Gamma$. Moreover, if
$\gamma$ is null-homologous and $S$ a Seifert surface, then
\begin{equation}\label{e:rcs}
\rot(\gamma)= \chi(S^+)-\chi(S^-),
\end{equation}
where $S^\pm$ is as in the definition of dividing set and $\chi(S^\pm)$ denotes the Euler characteristic.
\end{Pro}

If $\gamma$ is a Legendrian curve contained in a convex surface $S$, i.e.~not necessarily a boundary
component, $\gamma$ can be made to have a standard collar neighbourhood as depicted in Fig. \ref{f:collar}
(where $\gamma$ is a ruling curve in the interior); see \cite{kanda:1}.
Formula (\ref{e:twcs}) is also valid in this case.

Giroux pointed out that for convex surfaces, the dividing set, not
the particular characteristic foliation, essentially determines the contact structure in a neighbourhood.
Namely:
\begin{The}[Giroux's Flexibility Theorem, \cite{giroux:5}]\label{th:gflex} Consider a surface $S$, closed or
compact with Legendrian boundary, in a contact manifold $(M,\xi)$. Assume $\Gamma$ is a dividing set for the 
characteristic foliation $S_\xi$ and $\fol$ is another singular foliation on $S$ divided by $\Gamma$.
Then there is an isotopy $\phi: S\times[0,1]\to M$ of $S$ such that $\phi_0 = id$,
$\left(\phi_1(S)\right)_\xi = \phi_1(\fol)$, $\left(\phi_t(S)\right)_\xi$ is divided by $\Gamma$ for all
$t\in [0,1]$ and $\Gamma$ is fixed.
\end{The}

On the other hand, on a convex surface in a tight contact manifold, two dividing sets of a
characteristic foliation are isotopic. We will then, by slightly abusing language, refer to $\Gamma$ as
`the' dividing set.

As a consequence of Giroux's Flexibility Theorem, one can realize curves or arcs in a convex surface to be
Legendrian:
\begin{The}[Legendrian Realization, \cite{honda:1}] Consider a collection of disjoint properly embedded
closed curves and arcs $C$ on a convex surface $S$, which satisfies the following properties:
\begin{enumerate}
\item[(i)] $C$ is transverse to the dividing set $\Gamma$ of $S$ and every arc in $C$ begins and ends on
$\Gamma$,
\item[(ii)] every component of $S\setminus (\Gamma\cup C)$ has a boundary component which intersects
$\Gamma$,
\end{enumerate}
then there exists an isotopy $\phi:S\times [0,1]\to M$ such that $\phi_0 = id$, $\phi_t(S)$ are all convex,
$\phi_1(\Gamma) = \Gamma$ and $\phi_1(C)$ is Legendrian.
\end{The}

In a tight contact structure, the possibilities of dividing sets is rather restricted. Namely:
\begin{The}[Giroux's criterion, \cite{honda:1}]\label{th:gcrit} A convex surface (closed or compact with  Legendrian boundary)
$S$ other than the sphere $S^2$ has a tight neighbourhood if and only if no component of $S\setminus\Gamma_S$
bounds a disc. A convex sphere $S^2$ has a tight neighbourhood if and only if $\#\Gamma_{S^2} = 1$, i.e. if there
exists exactly one dividing curve.
\end{The}

\subsubsection{Edge-rounding}

Next, we describe how to smooth out two convex surfaces intersecting transversally along a common Legendrian
curve with a negative twisting number and moreover, to relate the dividing set of the two surfaces to the
dividing set on the smoothed surface.
\begin{Pro}[Edge-Rounding Lemma, \cite{honda:1}]\label{p:edger} Assume $S_1$ and $S_2$ are convex surfaces,
with convex collar boundary, intersecting transversely inside a contact manifold $(M,\xi)$ along a
common Legendrian boundary curve $\gamma$ with negative twisting number. Suppose that $S_1$ and $S_2$ are oriented such that smoothing the edge yields an oriented surface $S'$. Then the edge may be smoothed so that the dividing set of
$S'$ is obtained from the dividing sets on $S_1$ and $S_2$, the dividing curves connect such that positive
(negative) regions $S_1^\pm$ of $S_1$ connect to positive (negative) regions $S_2^\pm$ of $S_2$ as indicated in figure \ref{f:rounding}.
\end{Pro}
\begin{proof}
After possibly a small perturbation of a neighbourhood $N$ of $\gamma$ in $M$, we can consider the following
situation: Consider $\real^2\times (\real/\integer)$ with coordinates $(x,y,z)$, and contact 1-form $\alpha = 
\sin(2\pi n z)dx + \cos(2\pi n z)dy$ for some $n\in \integer^{>0}$. Locally, a neighbourhood of $\gamma$ is
contactomorphic to $N_\varepsilon= \{x^2+y^2\leq \varepsilon\}$ and $\gamma$ is given by $x=y=0$. The convex
surfaces become $S_1\cap N_\varepsilon = \{ x= 0, 0 \leq y <\varepsilon\} $ and $S_2\cap N_\varepsilon = \{ y= 0,
0\leq x < \varepsilon\}$  oriented by $\partial_x$ and $\partial_y$ respectively. The transverse vector field 
for $\{(x-\delta)^2 + (y-\delta)^2 = \delta^2\}\cap N_\delta$, $0<\delta<\varepsilon$, is the inward-pointing 
radial vector $-\partial_r$ for the circle $\{(x-\delta)^2 + (y-\delta)^2 = \delta^2\}$.
Take $S'= (S_1\cup S_2 \setminus N_\delta) \cup \{(x-\delta)^2 + (y-\delta)^2 = \delta^2\}\cap N_\delta$,
the dividing curve $z= \frac{k}{2n}$ on $S_1$ then connects to the dividing curve $z= \frac{k}{2n}-
\frac{1}{4n}$ on $S_2$, for $k= 0,\ldots,2n-1$; see Fig. \ref{f:rounding}.
\begin{figure}[tb]
\includegraphics[width=10cm]{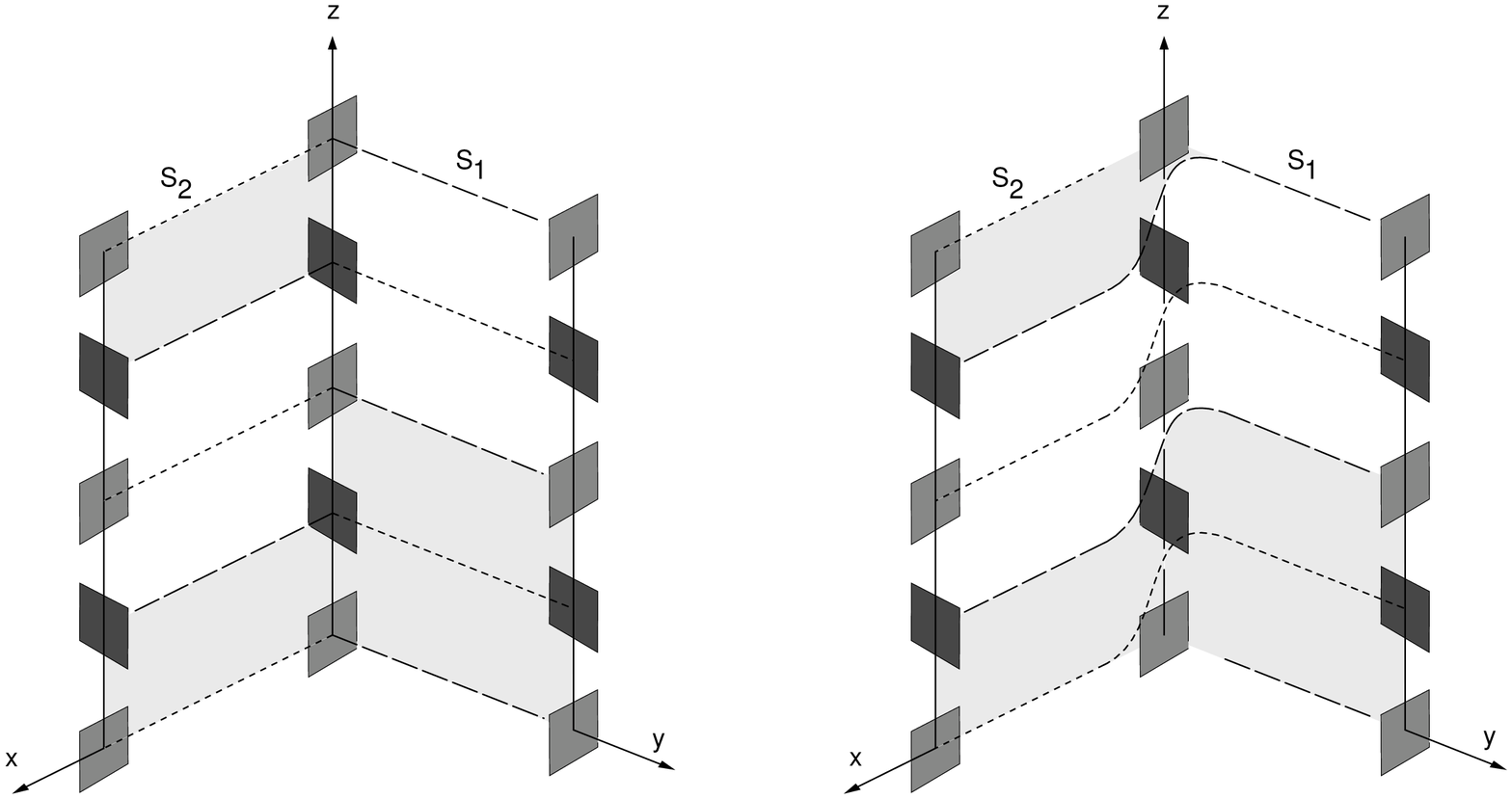}
\caption{Rounding edges. (Shaded regions are positive, dashed lines are dividing curves and Legendrian divides are dotted.)}
\label{f:rounding}
\end{figure}
\end{proof}

\subsubsection{Bypasses and alteration of the dividing set}

By perturbing a surface, one can alter the characteristic foliation to assume certain normal forms. 
On convex surfaces, we want to alter directly the dividing set. A crucial tool for this is the use of
{\it bypasses}, first exploited by Honda (see e.g. \cite{honda:1, honda:2}) with precursors in \cite{kanda:1}.

\begin{figure}[tb]
\includegraphics[width=7cm]{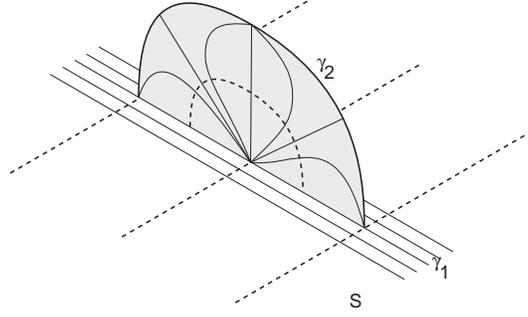}
\caption{A bypass. (Dashed lines are dividing curves on $S$.)}
\label{f:byp}
\end{figure}

\begin{Def} Assume $S \subset M$ is a convex surface. A {\it bypass} for $S$ is an oriented embedded disc
$D$ whose Legendrian boundary satisfies the following:
\begin{enumerate}
\item[(1)] $\partial D$ is the union of two arcs $\gamma_1$, $\gamma_2$ which intersect at their endpoints.
\item[(2)] $D$ intersects $S$ transversely along $\gamma_1$.
\item[(3)] along $\gamma_1$, there are three elliptic tangencies in the characteristic foliation $D_\xi$,
two of the same sign at the endpoints and one of different sign in the interior of $\gamma_1$
\item[(4)] along $\gamma_2$ there are at least three tangencies, all  have the same sign but alternating
indices.
\item[(5)] there are no interior singular points of $D_\xi$.
\end{enumerate}
See Fig. \ref{f:byp} for an illustration.
\end{Def}

Observe that all singular points on $\partial D$ have the same sign except the one elliptic point in the 
interior of $\gamma_1$. We call this the {\it sign} of the bypass.
The endpoints of $\gamma_1$ may be the same elliptic point, in this case we call $D$ a {\it degenerate}
bypass.

We first explain how to find bypasses and then give a discussion regarding how
bypasses are used to alter the dividing set of a convex surfaces. We discuss the cases used later in 
Section \ref{s:apl}, for a more complete discussion including applications, the reader may refer to the
literature; see \cite{honda:1,honda:2,honda:3,etnyre_honda:1, etnyre_honda:2,etnyre_honda:3}.

Assume $S$ is a convex surface with Legendrian boundary. After possibly perturbing $S$ we can further
assume that all boundary tangencies are half-elliptic (see Lemma 3.2 in \cite{honda:1}).
If $\tw(\gamma) = -n \leq 0$ for $\gamma \subset \partial S$, then the dividing curves
intersect $\gamma$ exactly $2n$ times. Suppose one of these dividing arcs is boundary-parallel, i.e. the arc
cuts off a half-disc which has no further intersections with $\Gamma_S$. A nearby arc in the complement,
parallel to this dividing arc, can be made Legendrian using the Realization Principle. After this, the arc bounds
a bypass. Thus, we have as a general principle:
\begin{Pro} Let $S$ be a compact surface having one Legendrian boundary with non-positive twisting number,
other than $D^2$ with $\tw(\partial D^2)=-1$. After possibly a small perturbation we can assume that $S$ is
convex and  all singular points of $S_\xi$ on $\partial S$ are half-elliptic. Suppose further $\gamma$ is a
boundary-parallel dividing curve. Then there exists a bypass which contains the half-disc cut off by $\gamma$. 
\end{Pro}

\begin{figure}[tb]
\includegraphics[width=10cm]{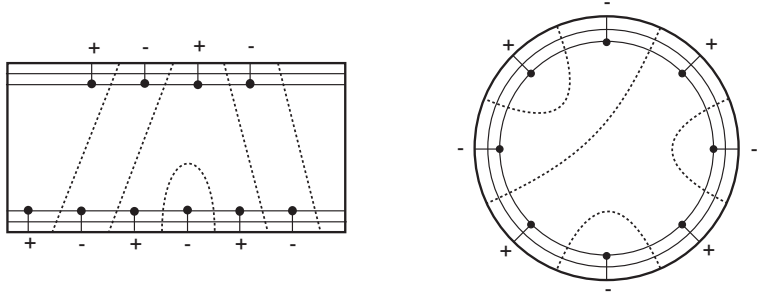}
\caption{A possible constellation of dividing curves on an annulus (left and right edge of the rectangle are identified) and a disc. (Dividing curves are dashed lines, on the right, we have $\tb(\partial D)=-4$)}
\label{f:div_together}
\end{figure}

In the sequel, we need bypass existence for two special surfaces: discs and annuli. We therefore consider these
special cases as discussed in \cite{honda:1}; see Fig. \ref{f:div_together} for examples.
\begin{Pro} Assume $D^2$ is a convex disc with Legendrian boundary lying inside a tight contact manifold and
$\tw(\partial D^2)=-n<0.$ After possibly a small perturbation we can assume that all tangencies at the boundary
are half-elliptic. Then $\Gamma_{D^2}$ consists of arcs which begin and end on $\partial D^2$. 
If $\tw(\partial D^2)<-1$, then there exists a bypass along $\partial D^2$.
\end{Pro}
\begin{Pro}[Imbalance Principle]\label{p:imbalance}
Assume $A=S^1\times[0,1]$ is a convex annulus with Legendrian boundary in a tight contact manifold.  After
possibly a small perturbation we can assume that all tangencies at the boundary are half-elliptic. If 
$\tw(S^1\times\{0\})< \tw(S^1\times\{1\}) \leq 0$, then there exists a bypass along $S^1\times \{0\}$.
\end{Pro}

\begin{figure}[tb]
\includegraphics[width=7cm]{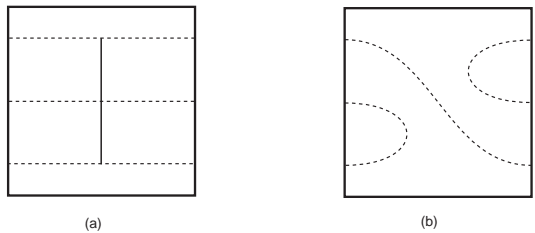}
\caption{Bypass attachment. Dividing curves on $Q$ before (a) and after (b) isotopy. The Legendrian line where the bypass is attached is drawn as thin vertical arc in (a).}
\label{f:b_attachment}
\end{figure}

Once we have a bypass for a convex surface, it can be used to manipulate the dividing set. The basic attachment
process is described as follows:
\begin{Pro}[Bypass attachment]\label{p:b_at} Assume $Q=[0,1]\times [0,1]$ is a convex square in a convex surface 
with three horizontal dividing arcs as in Fig. \ref{f:b_attachment} (a). If there exists a bypass $D$ for $Q$, along a vertical
Legendrian arc $\delta$, we can isotope $Q$ (fixing the boundary) by pushing $Q$ across $D$ such that the
characteristic foliation has a dividing set as shown in Fig. \ref{f:b_attachment} (b).
\end{Pro}

\begin{figure}[h]
\includegraphics[width=4cm]{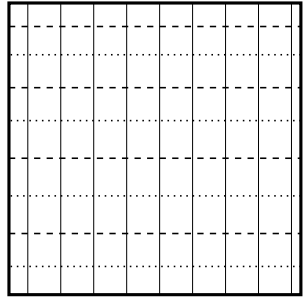}
\caption{A torus $T^2$ in standard form, in some identification with $\real^2/ \integer^2$, (Dividing curves are dashed, Legendrian divides are dotted horizontal and Legendrian  rulings are vertical solid lines) i.e. the sides are identified  and the bottom and top are identified.}
\label{f:stdT2}
\end{figure}

\begin{proof}
Because $Q$ is convex, we can consider an $I$-invariant one-sided neighbourhood $Q\times [0,\varepsilon]$,
for some $\varepsilon > 0$, such that $Q = Q\times \{\varepsilon\}$. Then, $A'=\delta\times [0,\varepsilon]$ is a
rectangle in standard form, (i.e. with horizontal linear characteristic foliation and parallel dividing arcs
in $[0,\varepsilon]$-direction), transverse to $Q$. Then, $A= A'\cup D$ is convex with piecewise smooth boundary.
The endpoints of the arc $\delta = A'\cap D$ are half-elliptic corners. In order to smooth the corners of $A$,
we convert half-elliptic points to full elliptic points. Finally, we apply the Pivot Lemma\footnote{
The Pivot Lemma allows to perturb a surface near an elliptic point such that any two transverse trajectories
become smooth; see \cite{eliashberg_fraser}.} to smooth the corners; see \cite{honda:1}. Because $A$ is convex,
we can take an $I$-invariant neighbourhood $N(A)= A\times[0,1]$. The boundary components $A_i= A\times\{i\}$,
$i=0,1$ are copies of $A$, i.e. have the same dividing set $\Gamma$. Both $A_0$ and $A_1$ are
oriented as boundary of $N(A)$ and therefore corresponding regions of $A_0\setminus\Gamma$ and
$A_1\setminus\Gamma$ have different signs. Now, using the Edge-Rounding Lemma \ref{p:edger}, we smooth out the
four edges of $N(A)\cup Q\times[0,\varepsilon]$ to obtain a surface $Q'$ with dividing set as in Fig.
\ref{f:b_attachment} (b), which completes the proof. See Fig. \ref{f:by_att} for an illustration.
\begin{figure}[tb]
\includegraphics[width=10cm]{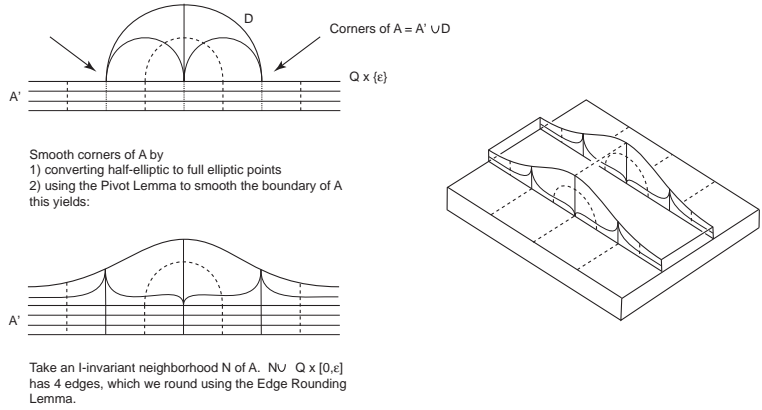}
\caption{Bypass attachment. An illustration of the proof. (Dashed lines are dividing curves.)}
\label{f:by_att}
\end{figure}
\end{proof}
In the sequel of this paper, we frequently encounter the situation where a bypass is attached along a torus.
We first describe a standard normal form for a convex torus in a contact structure and explain then the
consequences of the bypass attachment in this situation.  
On a convex torus $T^2$ in a tight contact manifold, we know by Giroux's criterion (Theorem \ref{th:gcrit}) that
no dividing curve bounds a disc. Therefore, the dividing set $\Gamma_{T^2}$ consists of $2n$ homotopic essential
parallel dividing curves and the number $n= \frac{1}{2}\#\Gamma_{T^2}$ is called the {\it torus division number}.
 Using some identification of $T^2$ with $\real^2/\integer^2$, the dividing curves have slope $s$, called the
{\it boundary slope} of the torus. Due to Giroux's Flexibility Theorem \ref{th:gflex}, we can deform the torus
$T^2$ inside a neighbourhood of $T^2\subset M$, fixing the dividing set $\Gamma_{T^2}$ so that the characteristic
foliation $T^2_\xi$ consists of a 1-parameter family of closed curves, called {\it Legendrian rulings}, of the
same slope $r$, called {\it ruling slope}. Each component of $T^2\setminus \Gamma$ contains a line of
singular points of slope $s$, called {\it Legendrian divide}. A convex torus in this (non-generic) form is
said to be in {\it standard form} (see Fig. \ref{f:stdT2} for an illustration).
An immediate consequence of Giroux's Flexibility Theorem is the following:
\begin{Pro}[flexibility of Legendrian rulings, \cite{honda:1}]\label{p:flexrul}
Assume $T^2$ is a convex torus in standard form, and, using coordinates in $\real^2/\integer^2$, has
boundary slope $s$ and ruling slope $r$. Then by a $C^0$-small
perturbation near the Legendrian divides, we can modify the ruling slope from $r\neq s$ to any other 
$r'\neq s$ ($\infty$ included).
\end{Pro}

\begin{figure}[tb]
\includegraphics[width=7cm]{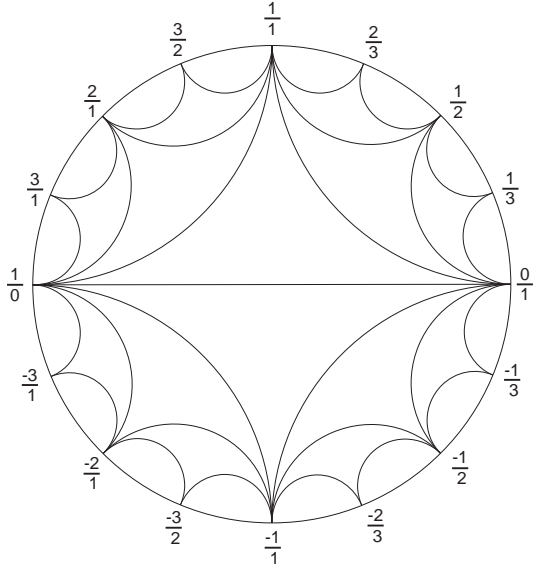}
\caption{The Farey tessellation of the hyperbolic disc.}
\label{f:farey}
\end{figure}

If a bypass is attached along some Legendrian ruling on $T^2$, we can push the torus across the bypass, using
the bypass attachment (Proposition \ref{p:b_at}), which yields a new torus with different boundary conditions.
If the torus division number $n$ of $T^2$ is greater than one, this will yield a torus with division number $n-1$.
In the case $n=1$ attaching a bypass does not change the torus division number but the boundary slope of the
torus; see \cite{honda:1}. In order to describe how the new boundary conditions are obtained from the old,
we first recall the Farey tessellation of the hyperbolic disc: Consider the hyperbolic unit disc
$\hyp = \{(x,y) \: : \: x^2 + y^2 \leq 1\}$. We label the point $(1,0)$ as $0=\frac{0}{1}$, the point $(-1,0)$
as $\infty = \frac{1}{0}$ and join the two points by an arc. Now label inductively points on
$S^1=\partial \hyp$ as follows (for $y>0$): assume we have labelled two points $\infty \geq \frac{p}{q}$,
$\frac{p'}{q'} \geq 0$ (where numerators and denominators are relatively prime). Label the point half
way between $\frac{p}{q}$ and $\frac{p'}{q'}$ along the shorter arc on $S^1$ by $\frac{(p+p')}{(q+q')}$. Connect two points
$\frac{p}{q}$ and $\frac{p'}{q'}$ by an arc if the corresponding shortest integral vectors form an integral
basis of $\integer^2$; see Fig. \ref{f:farey}.

\begin{The}[Honda, \cite{honda:1}]\label{t:byp_att}
Assume $T$ is a convex torus in standard form with $\#\Gamma=2$ and
boundary slope $s=s(T)$. If a bypass $D$ is attached to $T$ along a Legendrian ruling curve of slope
$r\neq s$, then the resulting convex torus $T'$ will have $\#\Gamma_{T'} = 2$ and boundary slope $s'$
which is obtained as follows: take the arc $[r,s] \subset \partial \hyp$ obtained by starting from $r$ and
moving counterclockwise until we hit $s$. On this arc, $s'$ is the point which is closest to $r$ and has
an edge from $s'$ to $s$.
\end{The} 

\subsection{Tight contact structures on basic blocks}

A key principle in the classification of tight contact structures on 3-manifolds is to cut along convex
surfaces to obtain simpler pieces on which the classification is known. In this subsection, we review the
basic properties of tight contact structures on various simple pieces referred to as {\it basic blocks}.

\subsubsection{The 3-ball}
The following key Theorem was proven by Eliashberg in \cite{eliashberg:1}:
\begin{The}\label{t:B3} Assume there exists a contact structure $\xi$ on a neighbourhood of $\partial B^3$
such that $\partial B^3$ is convex and $\#\Gamma_{\partial B^3} = 1$. Then there exists a unique extension
of $\xi$ to a tight contact structure on $B^3$, up to an isotopy relative to $\partial B^3$.
\end{The}

\subsubsection{The solid torus $S^1\times D^2$}\label{s:solidtorus}

Assume $\gamma\subset M$ is a Legendrian curve with a negative twisting number $\tw(\gamma)= n$ with respect
to some fixed framing. The {\it standard tubular neighbourhood} $N(\gamma)$ of $\gamma$ is defined as
solid torus $S^1\times D^2$ with coordinates $(z, (x,y))$  and contact 1-form $\alpha = \sin(2\pi n z)dx
+ \cos(2\pi n z)dy$ and $\gamma= \{(z,(x,y)) \: : \: x= y= 0\}$. With respect to the fixed framing of
$\gamma$, we may identify $\partial N(\gamma) = \real^2/\integer^2$ such that the meridian is $(1,0)^T$ and 
the longitude (fixed by the framing) is $(0,1)^T$. Then the boundary slope is $s(\partial N(\gamma)) = \frac{1}{n}$. 

In standard neighbourhoods of Legendrian curves the model standard tubular neighbourhood provides a unique
tight contact structure. This fact was used extensively by Kanda \cite{kanda:1}, and proved
(in a slightly different form) by Makar-Limanov in \cite{makar-limanov}; we refer to Theorem 6.7 in
\cite{etnyre:1}.
\begin{Pro} \label{unique} There exists a unique tight contact structure on $S^1\times D^2$ with a fixed convex boundary
with $\#\Gamma_{\partial(S^1\times D^2)} = 2$ and slope $s(\partial(S^1\times D^2))= \frac{1}{n}$, where $n$
is a negative integer. With the possibility of modifying the characteristic foliation on the boundary using the
Flexibility Theorem (Proposition \ref{p:flexrul}),
the tight contact structure is isotopic to the standard neighbourhood of a Legendrian curve with twisting number
$n$.
\end{Pro}

Decreasing the twisting number of a Legendrian curve is feasible, as commonly understood, by adding a `zigzag'
in the front projection; see \cite{etnyre:1}. Increasing the twisting number is not an easy task, but possible
in the presence of bypasses.

\begin{Pro}[Twist Number Lemma, \cite{honda:1}]\label{p:twi}
Consider a Legendrian curve $\gamma$ in a contact manifold $(M,\xi)$ with twisting number $n$ relative to a
fixed framing and $N$ a standard tubular neighbourhood of $\gamma$. If there exists a bypass attached to a Legendrian
ruling curve of $\partial N$ of slope $r$ and $\frac{1}{r}\geq n+1$, then there exists a Legendrian curve with
twisting number $n+1$ isotopic to $\gamma$. Notice that this isotopy cannot be a Legendrian isotopy 
because the twisting number changes.
\end{Pro}

Suppose we have given a solid torus $S^1\times D^2$ and an oriented identification of the boundary torus
$\partial(S^1\times D^2)$ with $\real^2/\integer^2$ so that $(1,0)^T$ corresponds to the meridian and $(0,1)^T$
corresponds to a longitude. Assuming the boundary of $S^1\times D^2$ is a torus in standard form with
torus division number one, the number of tight contact structures are determined by the boundary slope, i.e. the slope of
the dividing curves. More precisely:
\begin{The}[Theorem 2.3, \cite{honda:1}]\label{t:csoltor} Let $S^1\times D^2$ be a solid torus with convex boundary $T^2$ in standard form.
If $\#\Gamma_{T^2}=2$ and the boundary slope $s(T^2) = -\frac{p}{q}$, $p\geq q > 1$, $(p,q)=1$ and 
continued fraction expansion $$ -\frac{p}{q}= r_0 - \frac{1}{r_1 - \frac{1}{r_2 \cdots - \frac{1}{r_k}}},$$
with all $r_i < -1$, then there are
exactly $|(r_0+1)(r_1+1)\ldots(r_{k-1}+1)r_k|$ tight contact structures on $S^1\times D^2$ up to isotopy fixing
the boundary.
\end{The}

\begin{Pro}[Lemma 3.16, \cite{etnyre_honda:3}]\label{p:torins} Assume $S^1\times D^2$ has convex boundary with boundary slope $s<0$.
Then we can find a convex torus parallel to the boundary $\partial(S^1\times D^2)$ with any boundary
slope in $[s,0)$.
\end{Pro}

\subsubsection{Toric annuli $T^2\times [0,1]$}

Assume a toric annulus $T^2\times I$ is given in coordinates $(x,y,z)\in \real^2/\integer^2\times[0,1]$.
A standard tight contact structure is given by  $\alpha = \sin(\frac{\pi}{2} z)dx + \cos(\frac{\pi}{2} z) dy$. 
Note that the boundary $\partial (T^2\times I) = T_0 - T_1$ consists (after a perturbation) of 
convex tori with
boundary slope $0$ and $\infty$ respectively. It is not hard to see that the tori $T^2\times\{z\}$ are linearly
foliated and the boundary slopes decrease as $z$ increases. More generally, one obtains models for tight
contact structures on toric annuli with different slopes on the boundary $T_0$, $T_1$ by changing the chosen
interval $I$ on the $z$-axis.

\begin{Pro}[Proposition 4.16, \cite{honda:1}]\label{p:pendenze} Assume a toric annulus $T^2\times I$ has convex boundary in standard form and the boundary slope
on $T_i= T^2\times\{i\}$ is $s_i$, $i=0,1$ respectively. Then we can find convex tori parallel to $T_0$ 
with any boundary slope $s$ in $[s_1,s_0]$ (if $s_0<s_1$ this means $[s_1,\infty]\cup [-\infty,s_0]$).
\end{Pro}
On the other hand, consider a tight contact structure $\xi$ on a toric annulus $T^2\times I$ with convex
boundary and boundary slope $s_i= s(T_i)$, $i=0,1$. We say $\xi$ is {\it minimally twisting}
(in the $I$-direction) if every convex torus parallel to the boundary has slope $s \in [s_1,s_0]$.

We will outline the classification of tight contact structures in thickened tori $T^2 \times I$.
For a detailed description, we refer to Honda \cite{honda:1}. To state the Theorems, we first recall the
notion of the relative Euler class.
Consider a complex line bundle $\xi$ on a $3$-manifold $M$ with boundary 
$\partial M$.  Assume $\xi|_{\partial M}$ has a nowhere vanishing section $s$.
We may define the {\it relative  Euler class} $e(\xi, s) \in H^2(M, \partial M)$ as 
the obstruction to extending $s$ to the 
whole manifold. It is related to the Euler class by the following exact 
sequence: \\
$$\begin{matrix}
H^1(\partial M)& \stackrel{\Delta} \longrightarrow  & H^2(M, \partial M) & \to & H^2(M) & \to & H^2( \partial M) \\
      &                 & e(\xi ,s)   & \mapsto  & e(\xi) & \mapsto &  0
\end{matrix}$$

The following two Lemmas are useful for the calculation of the relative Euler
class of contact structures. The proofs are found in Section 4.2 of \cite{honda:1}.
\begin{Lem} \label{20}
Let $(M, \xi)$ be a contact manifold with convex boundary, and $s$ a fixed
section of $\xi |_{\partial M}$.
\begin{enumerate}
\item[(1)] If $\Sigma \subset M$ is a closed convex surface with positive (resp. 
negative) region $R_+$ (resp. $R_-$) divided by $\Gamma_{\Sigma}$, then 
$\langle e(\xi), \Sigma \rangle = \langle e(\xi ,s), \Sigma \rangle= 
\chi (R_+)- \chi (R_-)$.
\item[(2)] If $\Sigma \subset M$ is a compact convex surface with Legendrian 
boundary on $\partial M$ and regions $R_+$ and $R_-$, and $s$ is homotopic to 
$s'$ which coincides with $\dot{\gamma}$ for every oriented 
connected component $\gamma$ of $\partial \Sigma$, then $\langle e( \xi ,s), \Sigma 
\rangle = \chi (R_+)- \chi (R_-)$.
\end{enumerate}
\end{Lem}
\begin{Lem}
Let $(M, \xi)$ be a tight contact manifold with convex boundary consisting of tori. Then the relative Euler class $\langle e(\xi ,s), \Sigma \rangle$ is independent of the slope of the Legendrian rulings, where $s$ is a nonzero section of $\xi$ tangent to the ruling curves.
\end{Lem}
In the following we use the Euler class only in its relative form, Thus, when $\partial M$ is a union of convex tori in standard form, we will write $e(\xi)$ instead of $e(\xi ,s)$ with an abuse of notations if the section $s$ comes from a vector field tangent to the Legendrian rulings.

We say that two rational slopes in $\rational \cup \{ \infty \}$ are {\it consecutive} if they are
joined by an edge in the Farey tessellation. The very basic building blocks for contact structures are the
minimally twisting tight contact structures on $T^2 \times I$ whose boundary slopes $s_0$ and $s_1$ are
consecutive. Such contact structures are called {\it basic slices}. We have the following classification result
for basic slices.
\begin{The} [\cite{honda:1}, Section 4.3]
Given consecutive $s_0$ and $s_1$, there are, up to isotopy fixed on the boundary, two minimally twisting tight 
contact structures on $T^2 \times I$ with boundary in standard form, $\# \Gamma_{T_i}=2$ for
$i=0,1$, and boundary slopes $s_0$ and $s_1$.
The two contact structures are distinguished by their relative Euler class, and both can be
contact-embedded in a tight contact structure on $T^3$.
\end{The}

Let $v_0$ and $v_1$ be shortest integer vectors representing the slopes $s_0$ and
$s_1$ of a basic slice, such that $(v_1, v_0)$ is a positively oriented basis.
The possible relative Euler classes of the basic slices are the Poincar{\'e} duals of the homology
classes represented by $\pm (v_1 - v_0)$. For this reason, in what follows we will refer to the isotopy class
of a basic slice as its {\it sign}.
 
The basic slices are basic in the sense that any minimally twisting tight contact structure on $T^2 \times I$ can be decomposed into basic slices, as explained in the following Theorem:
\begin{The} [Lemma 4.12, \cite{honda:1}]
Given a minimally twisting, tight contact structure on $T^2 \times I$ with boundary slopes $s_0$ and $s_1$,
we can find a partition $0=t_0 < \ldots < t_k=1$ of $[0,1]$ such that $T^2 \times \{ t_i \}$ is a convex torus in standard form with slope $s_{t_i}$ for any $i=0,\ldots, k$. These slopes form a counterclockwise sequence in the arc $[s_1, s_0]$ with the property that $s_{t_i}$ and $s_{t_{i+1}}$ are consecutive and there is no edge in the Farey tessellation joining $s_{t_i}$ and $s_{t_{i+2}}$.

Moreover such a basic slices decomposition is minimal in the sense that any other basic slices decomposition of
the given tight contact structure on $T^2 \times I$ is a further decomposition of this one.
\end{The}
Note that the intermediate slopes of the basic slices decomposition depend only on $s_0$ and $s_1$ and are
independent of the isotopy class of the contact structure.

Conversely, given a decomposition into basic slices, we have the following gluing Theorem, which is a particular case of the more general Theorem 4.24 in
\cite{honda:1}.
\begin{The}  
Let $s_1<s_0$ be rational slopes. Then every choice of signs for the basic slices in the basic slices
decomposition associated to $s_0$ and $s_1$ realizes a minimally twisting, tight contact structure on
$T^2 \times I$ with slopes $s_0$ and $s_1$.
\end{The}
We observe that, unlike basic slices, in general these tight contact structures
cannot be contact-embedded into a tight contact structure on $T^3$. This is possible if and only if
the relative Euler class is $\pm (v_1-v_0)$.

We may ask when two different choices of signs for the basic slices give the 
same contact structure.
We say the basic slices in $T^2 \times [t_j, t_l]$ form a {\it continued fraction block} if there is
a slope $r$ such that there is an edge in the Farey tessellation joining $r$
and $t_i$ for all $j \leq i \leq l$.
To understand the origin of the name, see \cite{honda:1}, where this concept
appeared for the first time in this context. The importance of this notion 
comes from the fact that the sign of basic slices belonging to the same continued 
fraction block can be shuffled without affecting the isotopy type of the contact structure on $T \times I$.  
This is a nontrivial result whose proof can be found in \cite{honda:1}.

In this paper we will use the property of continued fraction blocks only in the
following case. Let $T^2 \times I$ carry a minimally twisting tight contact 
structure with boundary slopes $s_0 = - \frac{1}{n}$, for $n>0$, and $s_1 = \infty$.
Then all the basic slices of its decomposition belong to the same continued
fraction block, and therefore their signs can be shuffled.

As a result of the classification of basic slices, the basic slices 
decomposition, and this last fact about continued fraction blocks, we can 
derive the following classification Theorem for tight contact structures on 
$T^2 \times I$:
\begin{The} [Proposition 4.22, \cite{honda:1}]
The minimally twisting tight contact structures on $T^2 \times I$ with standard 
boundary,
$\# \Gamma_{T_i}=2$ and boundary slopes $s_0$ and $s_1$ are distinguished up to 
isotopy fixed on the boundary by their relative Euler class.  
\end{The} 

Actually, in  this Theorem Honda proves more than what he states: in fact he 
shows that, after normalising the boundary slopes to $s_0=-1$ and $s_1< -1$, the
tight contact structures are distinguished by the value their  relative 
Euler class takes at a horizontal annulus with Legendrian boundary. He proves
 this fact by showing that each basic slice in the decomposition gives a 
contribution to the value of the relative Euler class which is bigger than the 
sum of the contributions of the basic slices belonging to all the preceding
continued fraction blocks. The same arguments prove the following 
Corollary.
\begin{Cor} \label{utile}
The minimally twisting tight contact structures on $T^2 \times I$ with standard 
boundary, $\# \Gamma_{T_i}=2$ and boundary slopes $s_0=0$ and $s_1< -1$ are distinguished up to 
isotopy fixed on the boundary by the value their relative Euler class takes at a horizontal annulus with Legendrian boundary.
\end{Cor}

 
The number of minimally twisting tight contact structures on $T^2 \times I$ with 
fixed dividing set on the boundary is  finite, and is expressed  as a function of the continued fraction representation of $s_1$, after normalising $s_0$ to $-1$ by a change of coordinates, as in the case of solid tori in  Theorem \ref{t:csoltor}.

We mention also the case where the boundary slopes are equal. Then, either any convex intermediate torus $T \subset T^2 \times I$ has the same slope as the boundary tori, or there is a convex torus of slope $s$ for any $s \in \rational$. This is a consequence of Proposition \ref{p:pendenze}.

The tight contact structures on $T^2 \times I$ with boundary slopes $s_0=s_1$ such that
all the intermediate convex tori have the same slope are called {\em non-rotative}.

\section{Legendrian $(-1)$ surgery}

A very useful method to construct contact structures on a manifold is given by Legendrian surgery. Together with
the fact that, in particular, Legendrian $(-1)$ surgery preserves fillability one can construct tight contact
structures. Suppose $\gamma$ is a Legendrian curve in a contact manifold $(M, \xi)$ and we fix a framing $\fra$
of $\gamma$ such that the twisting number is zero $\tw(\gamma, \fra)=0$; for example take $\fra$ to be the
contact framing of $\gamma$. We then find a standard neighbourhood $N(\gamma)=S^1\times D^2$ with convex boundary so
that the dividing set $\Gamma_{\partial N(\gamma)}$ consists of two parallel curves. Take an oriented
identification $-\partial(M\setminus N(\gamma))=\partial N(\gamma)\simeq\real^2/\integer^2$ so that $(1,0)^T$
corresponds to the meridian and $(0,1)^T$ to the longitude given by a dividing curve. Thus the boundary slope
$s(\partial N(\gamma))$ is infinite and the meridian has slope zero. Let
$M'=(M\setminus N(\gamma))\cup_f N(\gamma)$, where $f : \partial N(\gamma) \to -\partial (M\setminus N(\gamma))$
is a diffeomorphism corresponding to
$$\left[\begin{array}{cc}  1 &  0  \\  -1 &  1  \end{array}\right] \in SL_2(\integer).$$
Topologically this corresponds to a $(-1)$ Dehn surgery along $\gamma$ with respect to the chosen framing, and
the contact structure can be glued together, after possibly adjusting the characteristic foliation, since the
dividing sets on $-\partial(M\setminus N(\gamma))$ and $f(\partial N(\gamma))$ are isotopic.
More generally one can define Legendrian $(r)$ surgery, for $r\in \rational$, as described in \cite{geiges_ding:1}.

Legendrian $(-1)$ surgery corresponds to a handle body construction in the sense of
\cite{gompf:1, weinstein} and thus preserves fillability. Recall that a contact manifold $(M,\xi)$ is called
{\it holomorphically fillable} if it is the oriented boundary of a compact Stein surface. For example the
standard tight contact structure on the three-sphere $S^3 \subset \complex^2$, given as the oriented plane field
of tangent complex lines, is holomorphically fillable. It is a remarkable result of Gromov \cite{gromov} and
Eliashberg \cite{eliashberg:2} that fillable contact structures are tight. Furthermore:
\begin{The}[Eliashberg, \cite{eliashberg:3}] If $(M',\xi')$ is obtained from a holomorphically fillable contact manifold $(M, \xi)$ by
Legendrian $(-1)$ surgery as described above, then $(M',\xi')$ is holomorphically fillable.
\end{The}

In fact both Theorems remain true for any notion of fillability (see \cite{etnyre_honda:2} for a survey), but
Legendrian $(-1)$ surgery does not preserve tightness in case of contact manifolds with boundary, as pointed out in
\cite{honda:3}.

\section{Applications}\label{s:apl}

In this section, we will show how to obtain the classification of tight contact structures in two specific cases:
the Brieskorn homology spheres $\pm\Sigma(2,3,11)$. Both manifolds are Seifert fibred spaces over
the sphere with three singular fibres.

Consider a Seifert manifold $M$ with three singular fibres over $S^2$. $M$ is described by Seifert
invariants $(\frac{\beta_1}{\alpha_1},\frac{\beta_2}{\alpha_2},\frac{\beta_3}{\alpha_3})$,
we refer to \cite{hatcher} for an introduction.

Assume $V_i$ are solid tori $S^1\times D^2$ with core curves $F_i$, $i= 1,2,3$. We identify
$\partial V_i$ with $\real^2/\integer^2$ by choosing $(1,0)^T$ as the meridional direction
and $(0,1)^T$ as a longitudinal direction.
Furthermore consider $S^1\times \Sigma$, where $\Sigma$ is a three-punctured sphere, i.e. a pair of pants. 
We identify each boundary component of $-\partial(S^1 \times \Sigma)$ with
$\real^2/\integer^2$ by setting $(0,1)^T$ as the direction of the $S^1-$fibre and $(1,0)^T$ as the direction
given by $-\partial(\{pt\}\times \Sigma)$. Then we obtain the Seifert manifold
$M(\frac{\beta_1}{\alpha_1},\frac{\beta_2}{\alpha_2},\frac{\beta_3}{\alpha_3})$
by attaching the solid tori $V_i$ to $S^1\times \Sigma$, where the attaching maps
$A_i : \partial V_i \to -\partial((S^1\times\Sigma)_i)$ are given by
$$A_i= \left[\begin{array}{cc} \alpha_i & \gamma_i \\ -\beta_i & \delta_i \\ \end{array}\right] \in
SL(2,\integer)$$

\begin{Rem} Note that we often refer to the same surface by different names. For example
$\partial V_i$ and $\partial(M\setminus V_i)$ denote the same torus, but the identification with
$\real^2/\integer^2$ is different.
\end{Rem}
\begin{Rem} Unless stated otherwise, properly embedded surfaces in a contact manifold are understood to
be convex, if possible.
\end{Rem}

Note that the three singular fibres $F_i$ may be isotoped to be Legendrian so that their twisting numbers  $n_i$
are particularly negative. Recall that a standard neighbourhood $V_i$ of $F_i$ with convex boundary has
boundary slope $\frac{1}{n_i}$. Furthermore we may assume that the ruling slope on $-\partial(M\setminus V_i)$
is infinite, thereby using the flexibility of Legendrian rulings (Proposition \ref{p:flexrul}).
Starting with this initial configuration, we do the following:

\begin{enumerate}
\item[(1)] In the first step, we try to maximise the twisting numbers of the singular fibres. For this, consider
a vertical annulus $A=S^1\times I$ with Legendrian boundary along ruling curves of two different tori $V_i$,
$V_j$. If the Imbalance Principle forces a bypass on $A$ we may apply the Twist Number Lemma
(Proposition \ref{p:twi}) to increase one of the twisting numbers $n_i$ or $n_j$. Repeating this process, two
different situations might occur.
        \begin{enumerate}
        \item[(a)] Either there exists a bypass on $A$ however we cannot apply the Twist Number Lemma. In this case, we
        can thicken the tori by attaching the bypass. In the cases below this yields an infinite boundary slope on a
        $-\partial(M\setminus V_i)$. Consider a vertical annulus from a Legendrian divide to the other two tori, we can
        thicken all three tori so that $s(-\partial(M\setminus V_i))=\infty$, $i=1,2,3.$
        \item[(b)] There exists no bypass on $A$. In this case a tubular neighbourhood of $V_i\cup V_j\cup A$ is a
        piecewise smooth torus with exactly four edges. Rounding the edges using the Edge-Rounding Lemma
        (Proposition \ref{p:edger}), we obtain a torus with boundary slope $s$, which can be thought of as the
        boundary of a neighbourhood of the third singular fibre $F_k$. In case $s < s(\partial(M\setminus V_k))$
        we can eventually increase the twisting number of $F_k$.
        \end{enumerate}
In either case, this process ends in a configuration with fixed boundary conditions on the basic blocks
$V_i$, $i= 1,2,3$ and $M\setminus\{V_1\cup V_2 \cup V_3\}$.
Combinations of tight contact structures on the basic blocks give a possible tight contact structure on $M$. 
Since there are finitely many tight contact structures on each basic piece, we obtain an upper bound on the
number of tight contact structures on $M$.
\item[(2)] In the second step we try to further analyse combinations of tight contact structures on the basic
blocks. Observe that if we were able to find further bypasses to thicken one $V_i$ such that $V_i$
contains a neighbourhood $V_i'$ of $F_i$ so that the boundary slope $s(\partial V_i')$ is zero, we would
find an overtwisted disc as meridional disc with boundary a Legendrian divide, and thus reduce the number of
potentially tight contact structures on $M$.

In the examples below, we are able to find further bypasses and eventually find an overtwisted disc in case
there exists a thickening so that $-\partial(M\setminus V_i)$ has infinite boundary slope.

\item[(3)]We finally construct tight contact structures by Legendrian surgery to show that the upper bound is sharp.
\end{enumerate}

\begin{Rem} Note that the strategy of constructing tight contact structures by Legendrian surgery yields
fillable contact structures. Hence this strategy is not successful in every case; see \cite{etnyre_honda:2} and
section \ref{s:conclusion}.
\end{Rem}

\subsection{The case $\Sigma(2,3,11)$}\label{s:example1}

In this subsection, $M$ denotes the Brieskorn homology sphere $\Sigma(2,3,11)$. This corresponds to the
manifold $M(\frac{1}{2},-\frac{1}{3},-\frac{2}{11})$ in the notation above. 

\begin{The}\label{th:main:1} On the Seifert manifold $\Sigma(2,3,11)$ there exist,
up to isotopy, exactly two tight contact structures, which are both
 holomorphically fillable.
\end{The}

The attaching maps are given by
$$
A_1 = \left[\begin{array}{cc}  2 &  1  \\   -1 &  0  \end{array}\right], \qquad
A_2 = \left[\begin{array}{cc}  3 & -1  \\    1 &  0  \end{array}\right], \qquad
A_3 = \left[\begin{array}{cc} 11 & -6  \\    2 & -1  \end{array}\right].
$$
Assume the singular fibres $F_i$ are (simultaneously) isotoped to Legendrian and further isotoped such that
their twisting numbers $n_i$ are particularly negative.
The standard neighbourhood of $F_i$ is denoted by $V_i$ and the slope of the dividing curves on $\partial V_i$ 
is $\frac{1}{n_i}$. Because $A_1\cdot(n_1,1)^T = (2n_1+1,-n_1)^T$,
$A_2\cdot(n_2,1)^T = (3n_2-1, n_2)^T$ and $A_3\cdot(n_3,1)^T = (11n_3-6, 2n_3-1)^T$, we calculate the boundary
slopes on $-\partial(M\setminus V_i)$ ($i= 1,2,3$) to be $\frac{-n_1}{2n_1 + 1}$, $\frac{n_2}{3n_2-1}$, and
$\frac{2n_3-1}{11n_3-6}$, respectively.

\subsubsection{Increasing twisting numbers of singular fibres}

We try to increase the twisting numbers of the singular fibres as far as possible.
As described above, we start by assuming the twisting numbers $n_i< 0$ are particularly negative.
\begin{Lem}\label{l:twinc} We can increase the twisting numbers $n_i$ of the singular fibres $F_i$, $i=1,2,3$,
up to $n_1=-1$, $n_2=n_3 = 0$.
\end{Lem}
\begin{proof} Using the flexibility of Legendrian rulings, we modify the Legendrian rulings on each
$\partial(M\setminus V_i)$ to have infinite slope.
Consider a vertical annulus $S^1\times I$ from $\partial(M\setminus V_1)$ to $\partial(M\setminus V_2)$ such 
that the boundary consists of Legendrian ruling curves on the tori. Observe that the boundary of this annulus
intersects the dividing curves on $\partial(M\setminus V_i)$ exactly $2(2n_1+1)$ and $2(3n_2-1)$ times respectively.

If $2n_1+1\neq 3n_2-1$, then, due to the Imbalance Principle (Proposition \ref{p:imbalance}), there exists
a bypass along a Legendrian ruling curve either on $\partial(M\setminus V_1)$ or $\partial(M\setminus V_2)$.
The Legendrian rulings on $\partial V_1$ have slope $-2$ and we can apply the Twist Number Lemma (Proposition \ref{p:twi})
to increase the twisting number of a singular fibre by one as long as $n_1 < -1$. A similar argument shows that we can use the
Twist Number Lemma to increase $n_2$ as long as $n_2 < 0$.

Assume $2n_1+1=3n_2-1$ and there exists no bypass on a vertical annulus $A=S^1\times I$ between
$\partial(M\setminus V_1)$ and $\partial(M\setminus V_2)$. We cut along the tori connected by $A$ and round the
corners using the Edge-Rounding Lemma: For this, observe that a neighbourhood of
$M\setminus (V_1\cup V_2\cup S^1\times I)$ is a piecewise smooth solid torus with four edges. Using the
Edge-Rounding Lemma (Proposition \ref{p:edger}), each rounding changes the slope by an amount
$-\frac{1}{4}\frac{1}{2n_1+1}$. Because there are four edges to round, we get on the boundary
$\partial(M\setminus V_1\cup V_2 \cup S^1\times I)$ the slope
$$\frac{-n_1}{2n_1+1} + \frac{\frac{2}{3}(n_1+1)}{2n_1+1} + \frac{-1}{2n_1+1} = -\frac{1}{3}\frac{n_1+1}{2n_1+1}.$$
Note that we identified this torus with $\real^2/\integer^2$ in the same way as
$\partial(M\setminus V_3)$. Since $A_3^{-1}\cdot(6n_1+3,n_1+1)^T = (3,-n_1+5)^T$, this corresponds to slope
$s=-\frac{1}{3}n_1+\frac{5}{3}$ when measured using $\partial V_3$.
Now $s>1$ for $n_1<0$ and we find a standard neighbourhood $V_3$ of $F_3$ with infinite boundary slope,
corresponding to $n_3 =0$. We remark that the boundary slope becomes $\frac{1}{6}$, when measured with respect to
$-\partial(M\setminus V_3)$.

Next, to increase the twisting number $n_2$, take a vertical annulus $S^1\times I$ from a
Legendrian ruling on $\partial(M\setminus V_2)$ to a Legendrian ruling on $\partial(M\setminus V_3)$. Observe
that if $n_2<-2$, we have $|3n_2-1|> 6$ and thus there exists a bypass on the $V_2$ side along 
a vertical Legendrian ruling, which allows us to increase $n_2$ up to $-1$ by the Twist Number Lemma.
A similar argument shows that we can increase $n_1$ up to $-2$.

Now the slopes on $-\partial(M\setminus V_i)$ are $-\frac{2}{3}$, $\frac{1}{4}$ and $\frac{1}{6}$ respectively.
Taking, once more, a vertical annulus between $V_1$ and $V_2$, we find a bypass due to the Imbalance Principle
and are finally able to increase $n_2$ to $0$ and $n_1$ to $-1$.
\end{proof}

We have now arrived at $n_1=-1$, $n_2 = n_3 = 0$. Note that the boundary slopes on $-\partial(M\setminus V_i)$ are
$-1$, $0$ and $\frac{1}{6}$ respectively. Take again a vertical annulus between $V_1$ and $V_2$. There are
two possibilities: Either there exists a bypass along both boundary components or not.
If there is a bypass, the cutting and rounding construction yields a torus of infinite slope. We use vertical
annuli from a Legendrian divide of this torus to thicken each $V_i$ to $V_i'$ s.t. $-\partial(M\setminus V_i')$
has infinite boundary slope.
In case there is no bypass, we perform a cutting and rounding construction on $V_1$ and $V_2$ as in the
proof of Lemma \ref{l:twinc} to obtain a further thickening of $V_3$  to $V_3'$ such that
$-\partial(M\setminus V_3')$ has boundary slope $0$.
We have shown that there are two possibilities, distinguished by whether or not there exists a thickening of
all $V_i$ such that the boundary slope with respect to $-\partial(M\setminus V_i)$ is infinite for $i=1,2,3$.
We will primarily be concerned with

\subsubsection{The case when a thickening to infinite slope exists}\label{s:ise}

We will now show that all possible tight contact structures arising in this case are overtwisted.
We do this by patching together meridional discs of two solid tori thus obtaining a surface with boundary on
the third torus and relate its dividing set to the dividing curves given on the discs. This may produce
a bypass which allows a further thickening, i.e. increasing the twisting and eventually becoming overtwisted.
In order to do this patching we have to examine the possible tight contact structures on the complement
of the singular fibres $S^1\times \Sigma \cong M\setminus(\cup_iV_i')$.
Consider $\Sigma = \{1\}\times \Sigma$. Each boundary component of $\Sigma$ intersects the dividing set of
the corresponding tori twice and therefore contains exactly two half-elliptic points. The following two Lemmata
are proven by Etnyre and Honda in \cite{etnyre_honda:1}. We enclose the proofs for the reader's convenience.

\begin{Lem}\label{l:conf}
        The dividing set on $\Sigma$ consists of arcs, each connecting two different boundary components.
\end{Lem}
\begin{proof}
Assume there is a boundary-parallel dividing arc as shown, for example, in cases (A) and (B) of
Fig. \ref{f:Sigma}.
\begin{figure}[tb]
\includegraphics[height=7cm]{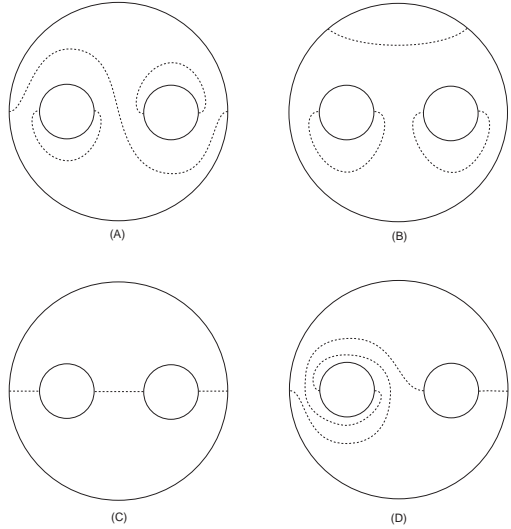}
\caption{Possible constellations of dividing curves on the pair of pants $\Sigma$. (Dividing curves are dashed lines.)}
\label{f:Sigma}
\end{figure}
This implies the existence of a bypass along some $\partial(M\setminus V_i')$. Attaching this bypass yields
a thickening $V_i''$ of $V_i'$ with slope $0$. Take a vertical annulus from a Legendrian dividing curve
on a $\partial(M\setminus V_j')$ ($j\neq i$) to $\partial(M\setminus V_i'')$. We find a bypass on this
annulus producing a further thickening to $V_i'''$ such that $\partial(M\setminus V_i''')$ has infinite
boundary slope. Therefore, by Proposition \ref{p:torins} we find a neighbourhood $V_i$ of $F_i$ so that the
boundary slope of $\partial V$ is zero. A meridional disc in $V$ whose boundary is a Legendrian divide on
$\partial V$ is an overtwisted disc.
Possible configurations of the dividing set on $\Sigma$ without boundary-parallel arcs are as shown 
in  Fig. \ref{f:Sigma} (C), up to twisting as shown in Fig. \ref{f:Sigma} (D).
\end{proof}
\begin{Lem} \label{30}
There exists a unique tight contact structure on $S^1\times \Sigma$, up to isotopy moving the boundary,
where the configuration of the dividing set on $\Sigma$ is given as in Lemma \ref{l:conf}.
\end{Lem}
\begin{proof}
We cut along $\Sigma$ and round the edges using the Edge-Rounding Lemma thus obtaining a solid two-handlebody. We can arrange the dividing set on the boundary so that two meridional discs intersect the dividing set
exactly twice; see Fig. \ref{f:Sigma_2}. Cutting along these two discs we obtain a three-ball.
Since there is a unique tight contact structure on the three-ball (Theorem \ref{t:B3}) and the dividing
curves on the surface we cut along are determined by the initial data, we must
 have a unique tight contact structure on $S^1\times \Sigma$.
\begin{figure}[tb]
\includegraphics[width=5cm]{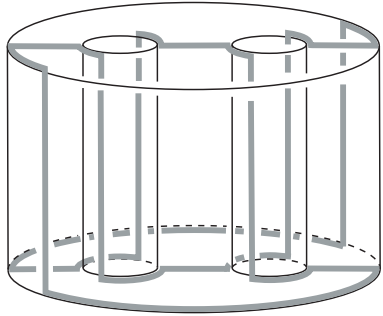}
\caption{The dividing set on the boundary of the two-handle. (The dividing curve, after rounding edges, is drawn as a thick curve.)}
\label{f:Sigma_2}
\end{figure}
\end{proof}

Using the flexibility of Legendrian rulings, we can choose the slopes on the $V_i'$ to be $0$, this is
possible since the boundary slopes, measured on $\partial V_i'$ are $-2$, $3$, and $\frac{11}{6}$, respectively. Take a
meridional disc $D_i$ in $V_i'$ such that $\partial D_i$ consists of a Legendrian ruling curve. On
$-\partial(M\setminus V_i')$, these discs have slopes $-\frac{1}{2}$, $\frac{1}{3}$ and $\frac{2}{11}$ respectively.

Since $\#(\partial D_i \cap \Gamma_{\partial(M\setminus V_i')})= 4$, $6$ and $22$, we obtain
$\tb(\partial D_i) = -2$, $-3$ and $-11$ and the possible constellations of dividing curves on $D_i$,
distinguished by their relative Euler number, i.e. the rotation number of $\partial D_i$ according to
$\rot (\partial D_i) = \chi(D_i^+) - \chi(D_i^-)$.

The two Lemmata above imply that the tight contact structure on $S^1\times \Sigma$ is contactomorphic to 
an $I$-invariant contact structure on a $T^2\times I$ with a neighbourhood of a vertical Legendrian curve of
zero twisting removed. View the $T^2\times I$ (minus $S^1\times D^2$) from the above Lemma as the region between
$\partial V_1'$ and $\partial V_2'$ (minus $\partial V_3'$) i.e. assume $T_0= \partial V_2'$ and
$T_1 = -\partial V_1'$. We write $T_t = T^2\times \{t\}$, $t\in [0,1]$.

Now pick three copies of meridional discs $D_1$ in $V_1$ and two copies of meridional discs $D_2$ in $V_2$.
Due to the $I$-invariance of $\xi$, we have a 1-parameter family of positive regions $(T_t)^+ = (T^2)^+
\times\{t\}$. Consider $D_1$ and $D_2$ such that $(D_2\cap T_0)^+ = \delta\times\{0\}$ and
$(D_1\cap T_1)^+ = \delta\times\{1\}$, where $\delta$ is a union of Legendrian arcs on $(T^2)^+$. Now 
$P= D_1\cup D_2\cup \delta\times [0,1]$ is a punctured torus. After smoothing the corners using the
Pivot-Lemma, $P$ has smooth boundary $\partial P\subset \partial V_3'$ with slope
$-\frac{1}{2}+\frac{1}{3}=-\frac{1}{6}$, measured using $\partial(M\setminus V_3')$.

\begin{figure}[tb]
\includegraphics[height=7cm]{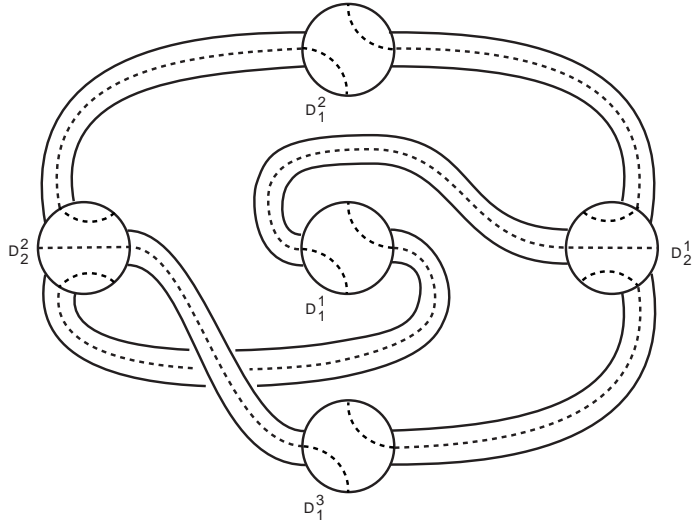}
\caption{Possible constellations of dividing curves on the punctured torus $P$, determined by the `signs' of the bypasses on the meridional discs. (Dividing curves are dashed lines.)}
\label{f:ptorus}
\end{figure}

{\it Case (1)} $D_1$ and $D_2$ have bypasses of the same sign. Then the dividing set on $P$ contains a
boundary-parallel curve, i.e. there exists a bypass on $P$. 
We can think of this bypass as attached on $-\partial(M\setminus V_3')$, along a ruling curve of slope
$\frac{1}{6}$. Attaching this bypass (Theorem \ref{t:byp_att}) yields a thickening of $V_3'$ to $V_3''$ such that
$-\partial(M\setminus V_3'')$ has boundary slope $1$.
Repeating the argument of Lemma \ref{l:twinc} shows that we can increase the twisting
numbers $n_1$ and $n_2$ to $-1$ and $0$ respectively. Thus we can thicken $V_3''$ further to $V_3'''$ such that
the boundary slope of $-\partial(M\setminus V_3''')$ is $0$. Recall that we started in a stage where the boundary slope
of $-\partial(M\setminus V_3)$ is zero. We assumed to find a thickening that this slope becomes infinite and showed that we
are then able to thicken further to obtain slope zero again. Thus we find a neighbourhood $V$ of $F_3$ so that $s(\partial V)=0$,
and a meridional disc in $V$ bounding a Legendrian divide is an overtwisted disc.

{\it Case (2)} $D_1$ and $D_2$ have different sign. Then there is a bypass on $D_3$ of the same sign as a
bypass on $D_1$ or $D_2$. Assume $D_1$ and $D_3$ contain bypasses of the same sign. A similar argument
as in Case (1) shows that patching eleven copies of $D_1$ and two copies of $D_3$ yields a surface whose
boundary is contained in $\partial(M\setminus V_2')$ with slope $-\frac{1}{2}+\frac{2}{11}=-\frac{7}{22}$. A bypass on each 
$D_1$ and $D_3$ joins to a bypass on the patched surface, if both have the same sign. Thus we find a bypass along
$V_2'$ and its attaching yields a thickening to $V_2''$ such that the boundary slope of
$-\partial(M\setminus V_2'')$ is $1$. Repeating the argument of Lemma \ref{l:twinc} shows that we can again
arrange $n_1=-1$ and $n_3=0$. A cutting and rounding construction gives a further thickening of $V_3'$ to
$V_3''$ such that $-\partial(M\setminus V_3'')$ has boundary slope $-1$ and hence contains an overtwisted disc.

Similarly, if $D_2$ and $D_3$ contain a bypass of the same sign, we patch together eleven copies of $D_2$ and
three copies of $D_3$ to obtain a surface whose boundary is contained in $\partial(M\setminus V_1')$ with
slope $\frac{1}{3}+\frac{2}{11}=\frac{17}{33}$.
We find a bypass and its attaching yields a thickening of $V_1$ to $V_1''$ such 
that $-\partial(M\setminus V_1'')$ has
boundary slope $0$. Make $n_2$ and $n_3$ again particularly negative and the same argument as in Lemma
\ref{l:twinc} shows that we can increase both $n_2$ and $n_3$ to $0$. Then, cutting and rounding along a
vertical annulus between $V_1$ and $V_2$ gives neighbourhood $V_3''$ of $F_3$ with boundary slope $1$
when measured using $-\partial(M\setminus V_3'')$. Since the boundary slope on $-\partial (M\setminus V_1'')$
is zero, we find by Proposition \ref{p:torins} a neighbourhood $V$ of $F_1$ so that $-\partial(M\setminus V)$
has infinite boundary slope. Take a vertical annulus $A=S^1\times I$ from a Legendrian divide on
$\partial(M\setminus V)$ to $\partial(M\setminus V_3'')$. There exists a bypass on $A$ whose attachment yields
a further thickening of $V_3''$ to $V_3'''$ where $\partial(M\setminus V_3''')$ has slope zero. Therefore, $V_3'''$
contains, by Proposition \ref{p:torins} a neighbourhood $V$ of $F_3$ so that the boundary slope of
$-\partial(M\setminus V)$ is $\frac {2}{11}$. A meridional disc in $V$ with boundary a Legendrian divide on
$\partial V$ is overtwisted.
Hence we have eliminated all possibilities in case there exists a thickening of the $V_i$ such that
$s(-\partial (M\setminus V_i))$ is infinite.

\subsubsection{The case when no thickening exists}

We are left now with the case when there exists no thickening of the standard neighbourhoods so that the
boundary slopes of the complements is infinite.
We have the following conditions: for the first singular fibre $F_1$ we obtained
twisting number $n_1 = -1$, hence a standard neighbourhood $V_1$ has boundary slope $-1$. Measured using 
$\partial(M\setminus V_1)$, the boundary slope is $1$, because $A_1 \cdot (-1,1)^T = (-1,1)^T$.
For the second singular fibre $F_2$ we obtained twisting number $n_2 = 0$ and hence a 
standard neighbourhood of $F_2$ has infinite boundary slope, which corresponds to slope $0$, when measured 
using $\partial(M\setminus V_i)$.
Lastly, for the third singular fibre, the twisting number is $n_3=0$ and the slope on
$-\partial (M \setminus V_3)$ is $\frac{1}{6}$. A cutting and
rounding construction along a vertical annulus between $V_1$ and $V_2$ yields a further thickening of $V_3$ such that
$-\partial(M \setminus V_3)$ has boundary slope $0$. 

In the first and second solid torus $V_1$ and $V_2$ there exists exactly one tight contact structure
as standard neighbourhood of Legendrian fibres. Because $A_3^{-1}\cdot(1,0)^T = (-1,-2)^T$ we find two tight
contact structures on $V_3$.

The remaining block is $\Sigma\times S^1 = M\setminus(V_1\cup V_2\cup V_3)$, where
we can arrange the boundary components of the pair of pants $\Sigma$ to be Legendrian along the boundary
components of $\Sigma\times S^1$. With this boundary conditions there exists exactly one tight contact
structure on this block. This is due to the following Lemma, as part of Lemma 5.1. in \cite{honda:2}
\begin{Lem} If for $S^1\times \Sigma$, where $\Sigma$ is a pair of pants, we have on the boundary tori
$\partial(S^1\times \Sigma) = T_1+T_2+T_3$ slopes $1$, $0$, $0$ respectively, then there exists
exactly one tight contact structure on $S^1\times \Sigma$ with no vertical Legendrian curve.
\end{Lem}

Thus there are at most two tight contact structures on $M$. In the next section, we use Legendrian 
surgery to see that there are two Stein fillable contact structures on $M$.

\subsubsection{Construction of a tight contact structure}

\begin{figure}[tb]
\includegraphics[height=7cm]{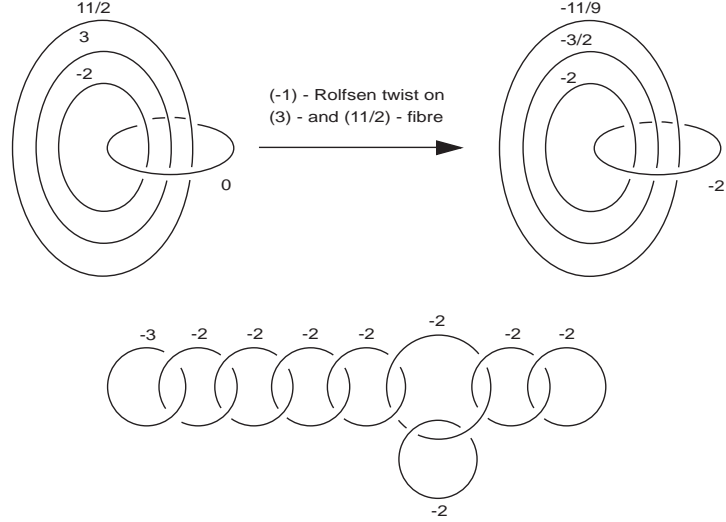}
\caption{Surgery description of the Seifert manifold. By Kirby calculus we obtain the description where all components of the link have coefficients $<-1$.}
\label{f:surgery}
\end{figure}

We will describe now how to establish a tight contact structure on $M$ by Legendrian surgery. 
The Seifert manifold $M(\frac{1}{2}, -\frac{1}{3},-\frac{2}{11})$ has a surgery description as the
left hand side of Fig. \ref{f:surgery}. By performing one (-1)-Rolfsen twist on the $(3)$ and $(\frac{11}{2})$ fibre,
we obtain the surgery description as shown on the right hand side; see \cite{gompf_stipsicz}.
Observe that we have the continued fraction expansions $-\frac{3}{2}=[-2,-2]$ and $-\frac{11}{9}=[-2,-2,-2,-2,-3]$, thus the
surgery description as at the bottom in Fig. \ref{f:surgery}. Because the surgery coefficients are $-2$ or $-3$ we may
conclude that there are exactly two Legendrian realizations of this link, where each component of the link will have
$\tb = -1$ and hence $r=0$, except one with $\tb =-2$ and $\rot = \pm 1$.
Note that the difference in the rotation number in a component of the two Legendrian realizations implies that the Chern
classes of the corresponding Stein surfaces are different, thus that the contact structures on the boundary
are non isotopic; see Theorem 4.20 in \cite{honda:1} and \cite{lisca_matic, lisca_matic:2} for details.

\subsection{The case $-\Sigma(2,3,11)$}

In this subsection $M$ will denote the Seifert manifold over $S^2$ with  three 
singular fibres with invariants $(- \frac{1}{2}, \frac{1}{3}, \frac{2}{11})$, corresponding to
the Brieskorn homology sphere $-\Sigma(2,3,11)$.

\begin{The}\label{th:main:2} On the Seifert manifold $-\Sigma(2,3,11)$ there exists,
up to isotopy, exactly one tight contact structure, which is holomorphically
 fillable.
\end{The}

Assume $V_i$ are tubular neighbourhoods of the singular fibres $F_i$, $i=1,2,3$,
and identify $M \setminus \cup V_i$ with $\Sigma \times S^1$, where $\Sigma$ is a pair of pants. We identify
$\partial V_i$ and $- \partial (M \setminus V_i)$ with $\real^2 / \integer^2$ as in the previous example.
 The gluing 
maps $A_i: \partial V_i \to - \partial (M \setminus V_i)$ are given by
$$
A_1= \left[ \begin{matrix} 2 & -1 \\  1 &  0  \end{matrix} \right],\qquad
A_2= \left[ \begin{matrix} 3 &  1 \\ -1 &  0  \end{matrix} \right],\qquad
A_3= \left[ \begin{matrix}11 &  6 \\ -2 & -1  \end{matrix} \right].
$$

\subsubsection{Increasing the twisting number of the singular fibres}
We begin by increasing the twisting number of the singular fibres as far as
possible in a similar way as in the previous example. We start by
assuming the singular fibres $F_i$ are simultaneously isotoped to Legendrian 
curves with twisting numbers $n_i$ very 
negative. The slopes of $\partial V_i$ are $\frac{1}{n_i}$, while the slopes of 
$-\partial (M \setminus V_i)$ are $\frac{n_1}{2n_1-1}$, $- \frac{n_2}{3n_2+1}$, and 
$- \frac{2n_3+1}{11n_3+6}$ respectively.

\begin{Lem}
We can increase the twisting numbers $n_1$ and $n_2$ up to $-2$, and the twisting
number $n_3$ up to $-1$.
\end{Lem}
\begin{proof}
Using Proposition \ref{p:flexrul}, we modify the Legendrian 
rulings on each $- \partial (M \setminus V_i)$ to have infinite slope. Consider a convex
annulus $A$ whose boundary consists of
Legendrian rulings of $- \partial (M \setminus V_1)$ and $- \partial (M \setminus V_2)$. If $2n_1-1 \neq 3n_2+1$ the 
Imbalance Principle provides a bypass along a Legendrian ruling either in 
$\partial (M \setminus V_1)$ or in $\partial (M  \setminus V_2)$. Using such a bypass we can apply the
Twist Number Lemma to increase the twisting number $n_i$ of a singular fibre by one as 
long as $n_1<0$ and $n_2< -1$.

If $2n_1-1=3n_2+1$, and there exist no bypasses on $A$, we get stuck in this 
operation. Suppose we are in this case: we cut along $A$ and round the edges,
obtaining a torus with slope $\frac{n_2}{6n_2+2}$ isotopic to
$\partial (M \setminus V_3)$. This slope corresponds to $- \frac{1}{2} n_2-2$ in $\partial V_3$, and is 
non-negative when $n_2 \leq -4$, therefore we can find a standard neighbourhood $V_3$ of $F_3$
with infinite boundary slope. This boundary slope becomes $- \frac{1}{6}$ if 
measured with respect to $- \partial (M \setminus V_3)$. To further increase $n_2$ take an 
annulus between $\partial(M \setminus V_2)$ and $\partial (M \setminus V_3)$.
If $n_2< -2$, we have $|3n_2+1|>6$ and thus there exists a bypass attached to
$\partial(M \setminus V_2)$ which allows us to increase $n_2$ by one so that
we can start again. In this way we can increase $n_1$ and $n_2$ up to $-2$.
When $n_1=n_2=-2$ the boundary slopes are $ \frac{2}{5}$ on $- \partial (M
\setminus V_1)$ and $- \frac{2}{5}$ on $- \partial (M \setminus V_2)$. By the Imbalance Principle, a convex 
annulus between $\partial (M \setminus V_2)$ and $\partial (M \setminus V_3)$ produces a bypass attached
to $\partial (M \setminus V_3)$ as long as $5 < |11n_3+6|$, thus we can use the Twist
Number Lemma to increase $n_3$ up to $-1$.
\end{proof} 

\begin{Lem} \label{1}
Let us suppose $n_1=n_2=-2$, $n_3=-1$ and $A$ is a convex vertical annulus
whose boundary consists of Legendrian rulings of $\partial (M \setminus V_1)$ 
and $\partial (M \setminus V_2)$. If $A$ has a boundary-parallel dividing 
curve, then the twisting numbers can be increased up to $n_1=0$, $n_2=n_3=-1$
and moreover there is a regular fibre with twisting number zero.
\end{Lem}
\begin{proof}
If there is a boundary-parallel dividing curve, then $A$ carries a bypass on 
each side after perturbing its characteristic foliation.
Using these bypasses,  we can further increase $n_1$ and $n_2$ up to $-1$. 
By the Imbalance Principle, we can find one more bypass in an annulus between $\partial (M \setminus V_1)$ and
$\partial (M \setminus V_2)$ on the side of $\partial (M \setminus V_1)$. This bypass increases the twisting 
number $n_1$ up to $0$. The slope of $- \partial (M \setminus V_1)$ is $0$, and the slope of
$- \partial (M \setminus V_2)$ is $- \frac{1}{2}$, therefore two possibilities for
an annulus $A$ between $\partial (M \setminus V_1)$ and $\partial (M \setminus V_2)$ are given: either  $A$ 
carries a bypass for $\partial (M \setminus V_1)$, or not. If such a bypass exists, then all the 
boundary slopes can be made infinite, and we can decrease the twisting $n_3$
to $-1$. If there is no such bypass, cutting along $A$ and rounding edges yields
a torus with slope $0$, which is $-2$ when measured in $\partial V_3$. In $V_3$
we find a convex torus with slope $- \frac{11}{6}$,
which corresponds to infinite slope in $- \partial (M \setminus V_3)$.  
\end{proof}
  
\subsubsection{The case when a thickening to infinite slope exists}
In this subsection we will show that there are no tight contact structures on 
$- \Sigma (2,3,11)$ with a regular fibre with twisting number zero. We suppose such 
a fibre exists and argue by contradiction. 

\begin{figure}
\includegraphics[height=5cm]{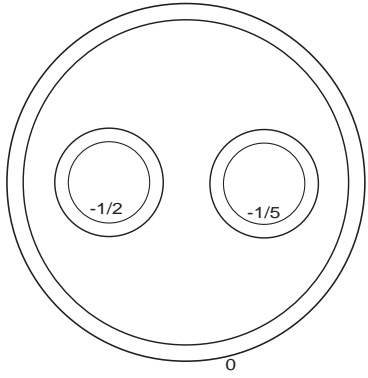}
\caption{$\Sigma_0$ with $\Sigma_1$ inside. The numbers refer to the slopes of the corresponding tori. The tori without indication of the slope are meant to have infinite slope.}
\label{pantalone.fig}
\end{figure}

Let $V_i $ be a standard neighbourhood of $F_i$. Then $M \setminus \cup V_i$ 
is diffeomorphic to $\Sigma_0 \times S^1$ for a pair of pants 
$\Sigma_0$, and has boundary slopes $0$, $- \frac{1}{2}$, and $- \frac{1}{5}$.
We use vertical annuli from a regular fibre with twisting number zero to 
thicken each $V_i$ to $V_i'$ such that $- \partial (M \setminus V_i')$ has infinite slope.  
We can find another pair of pants $\Sigma_1 \subset \Sigma_0$ such that 
$\Sigma_1 \times S^1$ is diffeomorphic to $M \setminus \cup V_1'$. The arguments in Lemma \ref{l:conf} 
apply to show that the dividing set of $\Sigma_1$ looks like in Figure \ref{f:Sigma} (C).
Note that $(\Sigma_0 \setminus \Sigma_1) \times S^1$ is the
 disjoint union of three thickened tori $T_i \times I$ such that $T_i \times
\{ 0 \} =- \partial (M \setminus V_i)$ and $T_i \times \{ 1 \} =- \partial (M
\setminus V_i')$ for $i=1,2,3$; see Fig. \ref{pantalone.fig}.
$T_1 \times I$ is a basic slice, $T_2 \times I$ is a union of
two basic slices $T_2 \times [0, \frac{1}{2}]$ with slopes $- \frac{1}{2}$ and
$-1$, and $T_2 \times [\frac{1}{2}, 1]$ with slopes $-1$ and $\infty$, and  
$T_3 \times I$ is a union of five basic slices $T_3 \times [\frac{i}{5}, \frac{i+1}{5}]$ 
for
$i=0, \ldots ,4$ with slopes $- \frac{1}{5-i}$ and $- \frac{1}{4-i}$. We  observe 
that the basic slices which 
compose the tight contact structures on the $T_i \times I$ belong to the same
continuous fraction block. 

The contact structures of the previous Lemma are described by the three numbers
$p_i$ of positive basic slices in $T_i \times I$, for $i=1,2,3$.

Our strategy to find overtwisted discs is cutting (a suitable sub-manifold of)
$M \setminus V_i \cup V_j$, for some $i$ and $j$, along a vertical annulus $A$ and rounding the 
edges of the cut open manifold to obtain a neighbourhood of the third singular
fibre with boundary slope zero. Slope zero on $V_i$ corresponds to $\frac{1}{2}$
on $- \partial (M \setminus V_1)$, to $- \frac{1}{3}$ on $- \partial (M \setminus V_2)$, and to $- \frac{2}{11}$
on $- \partial (M \setminus V_3)$. We will call these slopes {\em critical slopes}.
\begin{Rem}
This technique differs from the `disc patching' used in section \ref{s:ise} 
of the previous example.
We remark that both techniques work for each example and encourage the 
reader to perform both proofs with the technique not used.
\end{Rem}
The following Lemma allows us to control the dividing set of such an annulus $A$.
    
\begin{Lem} \label{modello}
Let $\Sigma$ be a pair of pants and $\xi$ a tight contact structure on $\Sigma
\times S^1$. Suppose that the boundary $- \partial (\Sigma \times S^1)=T_1 \cup T_2 \cup T_3$ consists of tori in 
standard form with $\# \Gamma_{T_i} =2$ for $i=0,1,2$, and slopes $s(T_0)= \frac{p_0}{q}$,
$s(T_1)= \infty$, $s(T_2)= \frac{p_2}{q}$. Suppose also that there exists a pair of 
pants $\Sigma' \subset \Sigma$ such that $\Sigma \times S^1$ decomposes as 
$\Sigma \times S^1 = \Sigma' \times S^1 \cup (T_0 \times I) \cup (T_2 \times I)$, where $\xi|_{\Sigma' \times S^1}$ is the tight 
contact structure with infinite 
boundary slopes described in Lemma \ref{30}, and $\xi|_{T_i \times I}$ is minimally 
twisting for $i=0,2$.

Suppose that one of the following holds:

\begin{enumerate}
\item $s(T_0)=s(T_1)=- \frac{1}{q}$ and $\xi|_{T_0 \times I}$ is isotopic to $\xi |_{T_2 \times I}$
\item $s(T_2)<0$ and $\xi|_{T_i \times I}$, for $i=0,2$, decomposes into basic slices of 
the same sign (i.e. with relative Euler class $\pm (q, p_i -1)$)
\end{enumerate}
Then there exists a convex annulus with Legendrian boundary consisting of 
vertical Legendrian rulings of $T_0$ and $T_2$ without boundary parallel dividing
curves.
\end{Lem}  
\begin{proof}
Let $\xi'$ be obtained in the following way from a tight contact structure on 
$T^2 \times I$ isotopic to 
$\xi|_{T_0 \times I}$. Remove a standard neighbourhood $U''$ of a vertical Legendrian 
ruling of a 
standard torus parallel to $T^2 \times \{ 0 \}$ and contained in its invariant 
neighbourhood, then thicken $U''$ to $U'$ with infinite boundary slope 
attaching the bypasses coming from the annulus between a vertical Legendrian 
ruling of $U''$ and a Legendrian divide of $T^2 \times \{ 1 \}$. By proposition 
\ref{p:pendenze}, there is a solid torus $U$ between $U''$ and $U'$ with 
boundary slope $\frac{p_2}{q}$ because $\frac{p_2}{q} \in [- \frac{1}{q}, - \infty)$. In
case (1) this operation is not necessary and we can simply take $U=U''$. In a 
similar way we can find a collar $C$ of $T^2 \times \{ 0 \}$ in $T^2 \times I \setminus U$ with 
boundary slopes $\frac{p_0}{q}$ and $\infty$.

We identify $T^2 \times I \setminus U$ to $\Sigma \times S^1$ so that $T^2 \times \{ 0 \}$ corresponds to $T_0$,   
$T^2 \times \{ 1 \}$ corresponds to $- T_1$, $\partial U$ corresponds to $T_2$, $C$ corresponds 
to $T_0 \times I$, and $U' \setminus U$
corresponds to $T_2 \times I$. A convex annulus with Legendrian boundary 
$A \subset (\Sigma \times S^1, \xi')$ between $T_0= T^2 \times \{ 0 \}$ and $T_2= \partial U$ contained in the 
invariant neighbourhood of $T^2 \times \{ 0 \}$ has no boundary parallel dividing 
arcs. To prove the Lemma we only need to show that $\xi$ and $\xi'$ are isotopic.
    
The dividing set of $\Sigma' \subset ( \Sigma \times S^1, \xi)$ cannot contain a boundary parallel 
dividing arc, otherwise there would be a convex torus with slope $0$ in $\Sigma' \times S^1$ (see 
the proof of Lemma \ref{l:conf}). By construction, $(\Sigma \times S^1, \xi)$ contact embeds in 
$T^2 \times S^1$ with a minimally twisting tight contact structure. A convex torus with slope $0$ around 
$T_2$ would give an overtwisted disc in $T^2 \times I$, while between $T_0 \times \{ 1 \}$ and 
$T_1$, both with infinite boundary slope, would contradict minimally twisting. 
This proves that $\xi |_{\Sigma' \times S^1}$ and $\xi' |_{\Sigma' \times S^1}$ induce the same dividing 
set on $\Sigma$, therefore they are isotopic by Lemma \ref{30}.

Take two vertical annuli with Legendrian boundary $A_i$, for $i=0,2$, between 
$T_i$ and $T_1$, we have 
$$\langle e(\xi'), A_0 \rangle = \langle e(\xi'), A_2 \rangle$$
because $[A_0]=[A_2]+[A]$ in $H_2(\Sigma \times S^1, \partial (\Sigma \times S^1))$ and $\langle e(\xi'), A \rangle=0$ by Lemma 
\ref{20}. Decompose $A_i =B_i \cup B_i'$, such that $B_i = A_i \cap T_i \times I$, and $B_i'=A_i \cap
\Sigma' \times S^1$ for $i=0,2$. We can suppose that all these annuli have Legendrian 
boundary. Since
$$\langle e(\xi' |_{\Sigma' \times S^1}), B_0' \rangle = \langle e(\xi' |_{\Sigma' \times S^1}), B_2' \rangle =0$$
we get 
$$\langle e(\xi' |_{T_0 \times I}), B_0 \rangle = \langle e(\xi' |_{T_2 \times I}), B_2 \rangle $$
and hence
$$ \langle e(\xi' |_{T_i \times I}), B_i \rangle = \langle e(\xi |_{T_i \times I}), B_i \rangle $$
for $i=0,2$ because $\langle e(\xi |_{T_2 \times I}), B_2 \rangle = \langle e(\xi |_{T_0 \times I}), B_0 \rangle = 
\langle e(\xi' |_{T_0 \times I}), B_0 \rangle$.
Applying Corollary \ref{utile} after a change of coordinates, we conclude
 $\xi |_{T_i \times I}$ is isotopic to $\xi' |_{T_i \times I}$ for $i=0,2$, and this ends the 
proof of the Lemma.
\end{proof}  
\begin{The}
 $M(- \frac{1}{2}, \frac{1}{3}, \frac{2}{11})$ carries no tight contact 
structure with a Legendrian regular fibre with twisting number zero.
\end{The}
\begin{proof}
For any choice of the number $p_i$ of positive basic slices in $T_i \times I$, for $i=0,1,2$, we will find a convex torus  
in $M \setminus \bigcup V_i$ with critical slope. Most of the possible tight contact structures
fall into one of the following cases.

{\it Case (1)}
We work between $V_2$ and $V_3$. If in $T_3 \times I$ there are two
basic slices with the same signs as the two basic slices in $T_2 \times I$, we can 
arrange them so that $T_3 \times [\frac{3}{5}, 1]$ is isotopic to $T_2 \times I$.
The manifold $M \setminus (V_1' \cup V_2 \cup V_3 \cup (T_3 \times [0,
\frac{3}{5}]))$ has boundary slopes $\infty$, $- \frac{1}{2}$, $- \frac{1}{2}$,
and by Lemma \ref{modello}, we
can find a convex annulus $A$ between $T_2 \times \{ 0 \}$ and $T_3 \times \{ \frac{3}{5} \}$
whose dividing curves go from a component of the boundary to the other one.
See Fig. \ref{case1.fig}.
After cutting along $A$ and rounding the edges, we obtain a torus with slope
$- \frac{1}{2}$ isotopic to $\partial (M \setminus V_1)$. Because the critical slope for 
$- \partial (M \setminus V_1)$ is $\frac{1}{2}$, $M$ is overtwisted. This case excludes all the
candidate tight contact structures except for the ones with $p_2=0$, $p_3=4,5$, 
or $p_2=1$, $p_3=0,5$, or $p_2=2$, $p_3=0,1$.

\begin{figure}[tb]
\includegraphics[height=5cm]{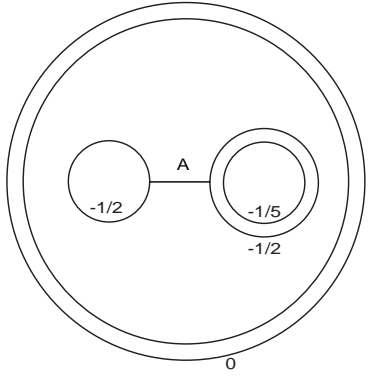}
\caption{Case 1. The tori without indication of the slope are meant to have infinite slope.}
\label{case1.fig}
\end{figure}

{\it Case (2)} We work between $V_1$ and $V_3$. Suppose that all the three basic
slices in $T_3 \times [\frac{2}{5} ,1]$ have the same sign as $T_1 \times I$.
We decrease the twisting number $n_1$ to $-1$ and take a standard 
neighbourhood $V_1''$ so that $- \partial (M \setminus V_1'')=T_1 \times \{ -1 \}$ has slope 
$\frac{1}{3}$. We can choose the sign of the basic slice 
$T_1 \times [-1,0]$  so that it is the same as $T_1 \times [0,1]$. To show this fact, 
embed $V_1$ in the 
standard $S^3$, and perform stabilisation there, see \cite{etnyre_honda:3}. 
We add $+1$ or $-1$ to the 
rotation number, according to the sign of the stabilisation, and it turns out 
that the 
relative Euler class of the contact structure on $V_1 \setminus V_2$ calculated
on a vertical annulus is the difference between the rotations.


The manifold $M \setminus (V_1'' \cup V_2' \cup V_3 \cup (T_3 \times [0, \frac{2}{5}]))$ has boundary 
slopes $\frac{1}{3}$, $- \frac{1}{3}$, 
$\infty$, and by Lemma \ref{modello} we can find a convex annulus $A$ between
$T_1 \times \{ -1 \}$ and $T_3 \times \{ \frac{2}{5} \}$ as in figure \ref{case2.fig} without 
boundary-parallel arcs, and by cutting along $A$
and rounding the edges we obtain a torus with slope $\frac{1}{3}$ isotopic to
$\partial (M \setminus V_2)$. This gives an overtwisted disk 
because $- \frac{1}{3}$ is the critical slope for $- \partial (M \setminus 
V_2)$. This case
excludes the contact structures with $p_1=0$ and $p_3 \leq 2$, or the contact 
structures with $p_1=1$ and $p_3 \geq 3$.

\begin{figure}[tb]
\includegraphics[height=5cm]{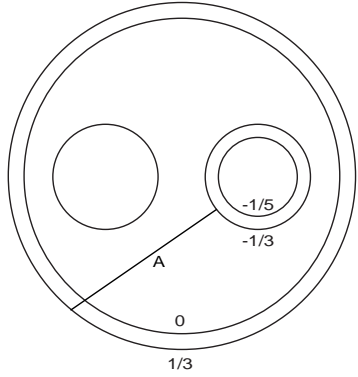}
\caption{Case 2. The tori without indication of the slope are meant to have infinite slope.}
\label{case2.fig}
\end{figure}

{\it Case (3)}
Now we work between $V_1$ and $V_2$. Suppose that the basic slices which 
compose $T_2 \times I$ have the same sign as $T_1 \times I$.
We decrease the twisting number of the singular
fibres $F_1$ and $F_2$  to $-2$ and take standard neighbourhoods $V_i''$, for
$i=1,2$, so that $- \partial (M \setminus V_1'')=T_1 \times \{ -1 \}$ has slope
$\frac{2}{5}$, $- \partial (M \setminus V_2'')=T_2 \times \{ -1 \}$ has slope 
$- \frac{2}{5}$, and all the basic slices in $T_1 \times [-1,1]$
and $T_2 \times [-1,1]$
have the same sign. The manifold $M \setminus (V_1'' \cup V_2'' \cup V_3')$ 
has boundary slopes $\frac{2}{5}$, $- \frac{2}{5}$, $\infty$, therefore 
by Lemma \ref{modello} we can find a convex 
 annulus $A \subset M \setminus V_3 \cup (T_3 \times I)$ between 
$T_1 \times \{ -1 \}$ and $T_2 \times \{ -1 \}$ without boundary-parallel
 dividing curves. See figure \ref{case3.fig}. Then, after cutting along $A$ 
and rounding the edges, we find a torus $V_3''$ 
isotopic to $\partial (M \setminus V_3)$ with  slope $\frac{1}{5}$. The 
thickened torus between $- \partial (M \setminus V_3)$ and $- \partial (M \setminus V_3'')$ has both boundary 
slopes $- \frac{1}{5}$, and contains 
 a torus with infinite slope, hence, by the classification Theorem for 
thickened tori, it contains an intermediate convex torus for each slope. In
particular, it contains a convex torus with slope $- \frac{2}{11}$, which is 
the critical slope for $- \partial (M \setminus V_3)$, therefore it gives an 
overtwisted disk around $F_3$. This case excludes the contact structures with
$p_1=0$, $p_2=0$ and $p_1=1$, $p_2=2$.

\begin{figure}[tb]
\includegraphics[height=5cm]{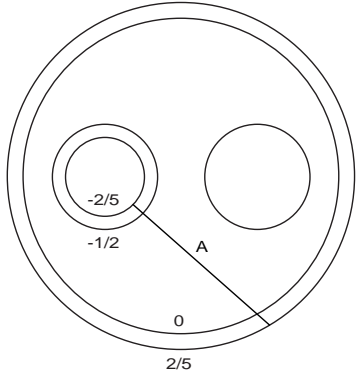}
\caption{Case 3. The tori without indication of the slope are meant to have infinite slope.}
\label{case3.fig}
\end{figure}

We can exclude the remaining candidate tight contact structures with $p_1=1$,
$p_2=1$, $p_3=0$, or $p_1=0$, $p_2=1$, $p_3=5$
in the following way. We go back to $V_2$ and $V_3$ and arrange the basic slices
in $T_2 \times I$ such that $T_2 \times [\frac{1}{2} ,1]$ has the same sign as the basic 
slices in $T_3 \times I$ and consider an annulus $A$ 
between $T_2 \times \{ \frac{1}{2} \}$ and $T_3 \times \{ \frac{4}{5} \}$ without boundary 
parallel dividing curves. Cutting along $A$ and rounding the edges we obtain a 
torus with slope $-1$ isotopic to $\partial (M \setminus V_1)$. This torus has slope $1$ in the
basis of $\partial V_1$, which  
 corresponds to increasing the twisting number of the singular fibre $F_1$ up
to $n_1=1$. Now we can  decrease $n_1$ again
choosing the sign of $T_1 \times I$, showing that the  contact 
structures with $p_1=1$, $p_2=1$, $p_3=0$, or $p_1=0$, $p_2=1$, $p_3=5$ are isotopic 
to the  contact structures
with $p_1=0$, $p_2=1$, $p_3=0$, or $p_1=1$, $p_2=1$, $p_3=5$ respectively, 
which have already been shown to be overtwisted in case (2).
\end{proof}
\subsubsection{The case when no thickening to infinite slope exists}
It remains to analyse only the case when the convex annulus $A$ between 
$\partial (M \setminus V_1)$ and $\partial (M \setminus V_2)$ of Lemma \ref{1} carries no bypasses. 
Since $- \partial (M \setminus V_1)$ has slope $\frac{2}{5}$ and $- \partial (M \setminus V_2)$ has slope 
$- \frac{2}{5}$, the dividing set of $A$ consists of $10$ dividing 
arcs going from one side of the annulus to the other. We say that an arc
in $A$ is horizontal if its algebraic intersection with 
$\Sigma \subset M \setminus \cup V_i \cong \Sigma \times S^1$ is zero.
\begin{Pro}
The manifold $M(- \frac{1}{2}, \frac{1}{3}, \frac{2}{11})$ carries at most
one tight contact structure.
\end{Pro}
\begin{proof}
Let $\phi_t$ be an isotopy of $A$ such that $\phi_0=id$ and $\phi_1(\Gamma_A)$ is a collection of 
horizontal arcs, and extend it to an isotopy on the whole $M$. We can consider
${\phi_1}_* \xi$ instead of $\xi$ and suppose without loss of generality that the 
dividing arcs of $A$ are horizontal. 

Cutting $M \setminus V_1 \cup V_2$ along $A$ and rounding the edges yields a solid torus 
with boundary slope $- \frac{1}{5}$, calculated with respect to the basis
of $- \partial (M \setminus V_3)$, which corresponds to slope $-1$ in the basis of $\partial V_3$ . By 
proposition \ref{unique}, there exists only one  tight contact 
structure up to isotopy on $V_3$ with this boundary condition, thus the 
candidate tight contact structure on $M(- \frac{1}{2}, \frac{1}{3}, 
\frac{2}{11})$ is unique up to isotopy. 
\end{proof}

\subsubsection{Construction of the tight contact structure}   
We will use Kirby calculus to show that this manifold can be represented as
Legendrian surgery on a link in $S^3$. This will prove that $M$ has at least 
one holomorphically fillable, and therefore tight, contact structure.
\begin{figure}[tb]
\includegraphics[width=10cm]{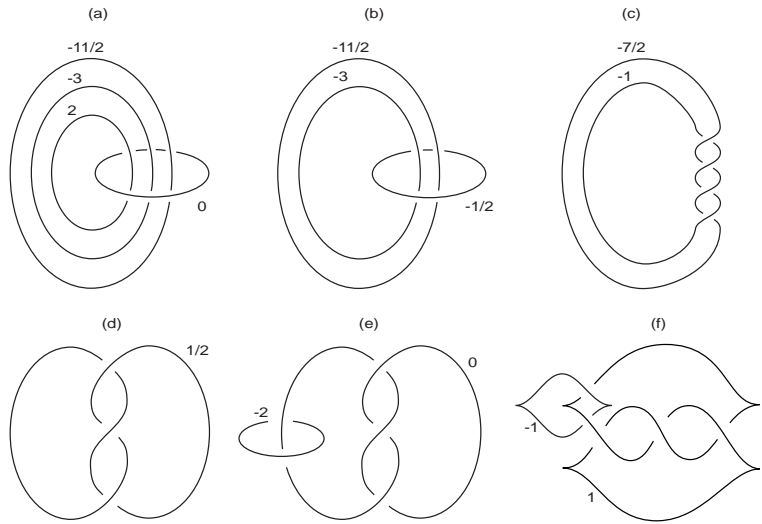}
\caption{Kirby calculus on $M$. The link as shown in (f) gives a Legendrian realization of the one in (e), the coefficients meaning the $\tb$ invariant of the respective component.}
\label{kirby1.fig}
\end{figure}

The structure of Seifert fibration of the manifold $M(- \frac{1}{2}, 
\frac{1}{3}, \frac{2}{11})$ gives the surgery presentation shown by (a) in
Fig. \ref{kirby1.fig}. By a slam-dunk on the $(2)$ component of 
the link, we obtain the link (b). Next, we perform a Rolfsen twist on the
$(- \frac{1}{2})$ component to obtain the link (c).
After another Rolfsen twist around the $(-1)$ component, we obtain diagram (d)
in Fig. \ref{kirby1.fig} and finally, after one inverse slam-dunk, we obtain
diagram (e). This link can be made Legendrian with the Thurston-Bennequin
invariant of each component one more than the surgery coefficient, see (f) for
 a Legendrian realization. 
This proves that $M(-\frac{1}{2}, \frac{1}{3}, \frac{2}{11})$ carries at least
one tight contact structure, and therefore ends the classification.

\section{What we can do, and what we can't do}\label{s:conclusion}

In this section we give a quick overview of similar results on the 
classification of tight contact structures on Seifert manifolds we are able to 
do with the techniques exposed in this paper, and point out some problems 
we have met. In addition to the Poincar{\'e} homology sphere with reversed 
orientation studied by Etnyre and Honda in \cite{etnyre_honda:1} and the two examples
treated in this paper, we are able to give a classification in several other cases.
Among them there are the manifold $M(\frac{1}{2}, - \frac{1}{2}, 
- \frac{1}{2})$, where the existence of a unique tight contact structure has 
been proven by S. Sch{\"o}nenberger; see \cite{schonenberger},
the Poincar{\'e} homology sphere $M(\frac{1}{2}, - \frac{1}{3}, - \frac{1}{5})$,
and the manifold $M(- \frac{1}{2}, \frac{1}{3}, \frac{3}{17})$. Moreover P.
Ghiggini in his Ph.D. thesis is working  towards the classification of tight contact
structures on Seifert manifolds over the torus.

We would like to point out that in the examples computed in this paper, like in
most of the other results mentioned above, the tight contact structures one 
obtains happen to be holomorphically fillable. In fact, the presence of 
possibly tight contact structures for which no Stein filling is known is a 
major source of difficulties in achieving a complete classification. The 
problems one has to face when dealing with non-holomorphically fillable
contact structures are twofold. 

The first problem is the proof of tightness, in fact we have very few 
techniques to do that, the main ones being  holomorphic or symplectic 
fillability and the gluing techniques developed by Colin \cite{colin:2} and Honda 
\cite{honda:3}. In particular, the gluing Theorems require the presence of 
incompressible surfaces, which do not exist in manifolds whose fundamental 
group is finite, such as small Seifert manifolds. An example is 
$M(- \frac{1}{2}, \frac{1}{3}, \frac{1}{4})$, which has been proved by Lisca to
be not fillable in any sense, and where only an upper bound for the number of 
tight contact structures is known.

The second problem is distinguishing the tight contact structures. For Stein 
fillable contact structures this is 
made easy by a result from 
Seiberg-Witten theory due to Lisca and Mati{\'c} \cite{lisca_matic} and Kronheimer and 
Mrowka \cite{kronheimer_mrowka}, which gives necessary conditions for holomorphically fillable
contact structures to be isotopic. Another way to distinguish contact 
structures is through their homotopy classification as $2$-planes fields
given by Gompf \cite{gompf:1} in terms of algebraic topological invariants.
An example of what happens outside the range of applicability of both methods 
is given by the family of manifolds $M(- \frac{1}{2}, \frac{1}{3}, 
\frac{2}{6n-1})$ for $n>3$, where we are able to give both an upper and a lower
bound on the number of tight contact structures, but the problem of 
determining whether some of them are isotopic or not remains open.

\section{Acknowledgement}

Both authors are very grateful to John Etnyre and Ko Honda for their encouragement,
steady support and interest. This work started during the Contact Geometry Quarter held at AIM and
Stanford University during Fall 2000. We would like to thank both institutions for providing support.

\bibliographystyle{amsplain}
\bibliography{contact}
\end{document}